\documentclass[11pt]{article}
\usepackage{amssymb,amsmath,booktabs,hyperref,mathtools}
\usepackage{graphics,authblk}
\usepackage{xcolor}
\usepackage{longtable,nicefrac}
\usepackage{graphicx}
\usepackage{subcaption}
\usepackage{multirow}
\usepackage{tabularx,enumerate,cases}
\usepackage{rotating}
\usepackage{dcolumn,accents}
\usepackage{pdflscape}
\usepackage{array}
\usepackage{lscape}
\usepackage{anysize}
\usepackage{color}
\usepackage{amsthm}
\usepackage{listings}
\newcommand{\mathleft}{\@fleqntrue\@mathmargin0pt}
\newcommand{\mathcenter}{\@fleqnfalse}

\providecommand{\keywords}[1]{\textbf{\textit{Keywords: }} #1}
\usepackage{float}
\usepackage{soulutf8}
\usepackage{adjustbox}
\usepackage{colortbl}
\usepackage{blindtext,nameref}

\newtheorem{theorem}{Theorem}
\newtheorem{corollary}{Corollary}
\newtheorem{proposition}{Proposition}
\newtheorem{lemma}{Lemma}
\theoremstyle{definition}
\newtheorem{definition}{Definition}
\theoremstyle{remark}
\newtheorem*{remark}{Remark}
\graphicspath{{photos/}}
\allowdisplaybreaks
\begin{document}
	\title{An Analytical Framework for the Linear Best-Worst Method and its Application to Achieve Sustainable Development Goals--Oriented Agri-Food Supply Chains}
	\author{Harshit M. Ratandhara, Mohit Kumar}
	\date{}
	\affil{Department of Basic Sciences,\\ Institute of Infrastructure, Technology, Research And Management, Ahmedabad, Gujarat-380026, India\\ Email: harshitratandhara1999@gmail.com, mohitkumar@iitram.ac.in}
	\maketitle
	\begin{abstract}
		The Best-Worst Method (BWM) has emerged as a prominent multi-criteria decision-making method for determining the weights of the decision criteria. Among various BWM models, this research focuses on the linear model of the BWM. This model calculates weights by solving an optimization problem, necessitating optimization software. In this article, we present a novel framework that solves this optimization model mathematically, yielding an analytical expression for the resultant weights, thus eliminating the requirement for an optimization software. The proposed approach enhances both the conceptual clarity of the underlying optimization process and the computational efficiency of the model. Based of this framework, we demonstrate the model's limited response to data variations, i.e., its lower data sensitivity. We also compute the values of consistency index for the linear BWM, which are required to calculate the consistency ratio - a consistency indicator used for assessing inconsistency in input data. Finally, we illustrate the validity and applicability of the proposed approach through five numerical examples and a real-world case study that ranks eighteen drivers across three categories - Industry 4.0, sustainability, and circular economy - in relation to sustainable development goals-driven agri-food supply chains.
	\end{abstract}
	\section{Introduction}
	Multi-Criteria Decision-Making (MCDM) is one of the key branches of operations research that aids decision-makers in evaluating complex problems, particularly those involving conflicting criteria. A critical step in addressing any decision-making situation is determining the weights of the decision criteria. Some of the most salient MCDM methods for deriving criteria weights include the Analytic Hierarchy Process (AHP) \cite{saaty1994make}, the Analytic Network Process (ANP) \cite{saaty2004decision}, and the Best-Worst Method (BWM) \cite{rezaei2015best}.\\\\
	The BWM is a distance-based method that utilizes pairwise comparisons for determining the criteria weights \cite{rezaei2015best}. The original model of BWM employs maximum deviation as a distance function and derives weights by solving a nonlinear minimization problem, thus known as the nonlinear BWM \cite{rezaei2015best}. Subsequent developments have introduced alternative BWM models, including the multiplicative BWM \cite{brunelli2019multiplicative}, the Euclidean BWM \cite{kocak2018euclidean}, and the goal programming-based BWM \cite{amiri2020goal}, each involving different distance functions. The BWM has also been extended to group decision-making scenarios \cite{safarzadeh2018group}. Some other theoretical developments regarding BWM are as follows. Rezaei \cite{rezaei2016best} applied interval analysis to rank interval-weights in the nonlinear BWM. He also introduced the concentration ratio to measure the dispersion of these interval-weights \cite{rezaei2020concentration}. The estimation of the accuracy of the decision data is a fundamental requirement for any MCDM method. The BWM framework incorporates consistency estimation through the Consistency Ratio (CR), calculated using the Consistency Index (CI) and the optimal objective value from the associated optimization problem. So, CR is an output-based consistency indicator that assesses the accuracy of the given comparisons only after the completion of full computation process \cite{rezaei2015best}. To address this limitation, Liang et al. \cite{liang2020consistency} proposed an input-based CR capable of providing real-time feedback to the decision-makers. They further established threshold values for this input-based CR to determine permissible inconsistency levels in pairwise comparisons. Mohammadi and Rezaei \cite{mohammadi2020bayesian} introduced the Bayesian BWM to examine the BWM framework probabilistically and derive aggregated criteria weights across multiple DMs. Lei et al. \cite{lei2022preference} addressed the challenge of improving the consistency of given preferences. They developed an optimization model that suggests optimal preference modifications, ensuring both ordinal consistency and an acceptable level of cardinal consistency. For deriving priority weights in BWM, Xu and Wang \cite{xu2024some} presented eleven methods for individual DMs and nine for group decision-making scenarios. Additionally, they discussed six inconsistency indices to assess the consistency of pairwise comparisons. Corrente et al. \cite{corrente2024better} introduced the parsimonious BWM, an extension of the nonlinear BWM, designed to facilitate priority determination for alternatives in situations involving a large number of alternatives. The BWM has also been combined with various other MCDM methods such as BWM-TOPSIS \cite{youssef2020integrated}, BWM-VIKOR \cite{dawood2023novel}, BWM-SERVQUAL \cite{rezaei2018quality}, and the Best-Worst Tradeoff (BWT) method \cite{liang2022best}. Additionally, the method has been extended to generalizations of classical sets, such as fuzzy sets \cite{guo2017fuzzy,ratandhara2024alpha}, intuitionistic fuzzy sets \cite{wan2021novel, mou2016intuitionistic}, and hesitant fuzzy sets \cite{ali2019hesitant}, to account for uncertainty in decision data.\\\\
	To address the non-uniqueness of resultant weight sets in the nonlinear BWM, Rezaei \cite{rezaei2016best} proposed a novel model of BWM following a similar philosophy to the nonlinear BWM. In this model, weights are derived by solving a linear minimization problem, hence it is known as the linear BWM. Due to its ability to produce a unique weight set, the linear BWM has become one of the most widely used BWM models in real-world scenarios, including sustainability assessment of supply chain \cite{ahmadi2017assessing}, evaluation of eco-industrial park \cite{zhao2018comprehensive}, ranking of barriers to energy efficiency \cite{gupta2017developing}, and location selection \cite{kheybari2020sustainable}. Despite the model's wide applicability, two key research gaps remain: the reliance on optimization software for solving the minimization problem reduces its time efficiency, and the lack of a systematic approach for the calculation of CI prevents its native consistency estimation capability.\\\\
	Wu et al. \cite{wu2023analytical} and Ratandhara and Kumar \cite{ratandhara2024analytical} analytically solved the optimization problems of the nonlinear BWM and the multiplicative BWM respectively, eliminating the requirement for optimization software in these models. Following a similar approach, in this article, we introduce a novel framework to derive an analytical expression for optimal weights in the linear BWM. To achieve this, we first partition the feasible region of the optimization problem for the linear BWM, i.e., the collection of all normalized weight sets, on the basis of the relationship between the weight assigned to the worst criterion and the objective function value of the weight set. We then establish some key properties of these block sets, including non-emptiness and relationships between weights within a specific block, which lead to a closed-form expression for the optimal weights. Beyond improving computational efficiency, this expression reveals the model's low sensitivity to data variations, exposing a key limitation of the model. We also compute CI using this expression. We then validate the proposed approach through five numerical examples. Furthermore, we demonstrate its applicability using a real-world case study involving the ranking of Industry 4.0, sustainability, and circular economy factors within sustainable development goals-driven agri-food supply chains.\\\\
	The rest of the article is organized as follows: Section 2 introduces fundamental concepts and provides a brief overview of the linear BWM. Section 3 develops an analytical framework for linear BWM, including derivation of an analytical expression for optimal weights, sensitivity analysis based on thi expression, calculation of CI, and validation of the approach through five numerical examples. Section 4 demonstrates a real-world application of the proposed methodology. Finally, Section 5 presents concluding remarks and discusses potential future research directions.
	\section{Basic Concepts and Linear Best-Worst Method in Brief}
	To resolve the issue of multiple resultant weight sets in the nonlinear BWM, Rezaei \cite{rezaei2016best} developed the linear BWM, which follows the same fundamental philosophy as the nonlinear BWM while ensuring unique resultant weights. This section covers some key concepts related to the linear BWM and provides a concise overview of the model.
	\subsection{Preliminaries}
		Let $C=\{c_1,c_2,\ldots,c_n\}$ be the set of decision criteria, and let $D=\{c_1,c_2,\ldots,c_n\}\setminus\{c_{b},c_{w}\}$ throughout the article. Whenever there is no ambiguity, we use the abbreviated notations $C=\{1,2,\ldots,n\}$ and $D=\{1,2,\ldots,n\}\setminus\{b,w\}$.\\\\
		The Pairwise Comparison System (PCS) is the pair $(A_b,A_w)$, where $A_b=(a_{b1},a_{b2},\ldots,a_{bn})$ is the best-to-other vector and $A_w=(a_{1w},a_{2w},\ldots,a_{nw})^T$ is the other-to-worst vector. Here, $a_{ij}$ represents the relative preference of the $i^{th}$ criterion over the $j^{th}$ criterion.
	\begin{definition} \cite{rezaei2015best}
		A PCS $(A_b,A_w)$ is said to be consistent if $a_{bi}\times a_{iw}=a_{bw}$ for all $i\in D$.
	\end{definition}
	\begin{theorem}\cite{wu2023analytical}\label{4accurate}
		The system of linear equations
		\begin{equation}\label{4system}
			\begin{split}
			&\frac{w_b}{w_i}=a_{bi},\quad \frac{w_i}{w_w}=a_{iw},\quad \frac{w_b}{w_w}=a_{bw},\ i\in D,\\
			&w_1+w_2+\ldots+w_n=1
			\end{split}
		\end{equation}
		has a solution if and only if $(A_b,A_w)$ is consistent. Furthermore, if solution exists, then it is unique and is given by
		\begin{equation}\label{accurate_weights}
			w_i=\frac{a_{iw}}{\displaystyle\sum_{j\in C}a_{jw}}=\frac{1}{a_{bi}\displaystyle\sum_{j\in C}\frac{1}{a_{bj}}},\ i\in C.
		\end{equation}
	\end{theorem}
	\subsection{Linear BWM}
	In the linear BWM, the decision-making process begins with the selection of the best criterion $c_b$ and the worst criterion $c_w$ from $C$. The decision-maker then constructs $(A_b,A_w)$ by providing the relative preferences of the best criterion over each criterion and the relative preferences of each criterion over the worst criterion. These preferences are usually provided as linguistic terms, which are then quantified using some established scale, such as Saaty scale \cite{saaty1994make}. Based on this PCS, the following minimax problem is formulated \cite{rezaei2016best}.
	\begin{equation}\label{4optimization_1}
		\begin{split}
			&\min\max\{|w_b-a_{bi}\times w_i|, |w_i-a_{iw}\times w_w|,|w_b-a_{bw}\times w_w|: i\in D\}\\
			&\text{subject to: }  w_1+w_2+\ldots+w_n=1 \text{ and } w_j\geq 0 \text{ for all }j\in C.
		\end{split}
	\end{equation}
	An optimal solution of this problem gives an optimal weight set. Now, consider the following minimization problem.
	\begin{equation}\label{4optimization_2}
		\begin{split}
			&\min \epsilon\\
			&\text{subject to: }\\
			&|w_b-a_{bi}\times w_i|\leq \epsilon,\ |w_i-a_{iw}\times w_w|\leq \epsilon,\ |w_b-a_{bw}\times w_w|\leq \epsilon,\\
			&w_1+w_2+\ldots+w_n=1 \text{ and } w_j\geq0 \text{ for all }i\in D \text{ and } j\in C.
		\end{split}
	\end{equation}
	Problem \eqref{4optimization_2} is an equivalent formulation of problem \eqref{4optimization_1}. This linear problem involves $n+1$ variables $w_1,w_2,\ldots,w_n$ and $\epsilon$. The optimal values $w_1^*,w_2^*,\ldots,w_n^*$ form an optimal weight set, while $\epsilon^*$, which is also the optimal objective value, estimates the accuracy of this weight set.
	\section{An Analytical Framework for the Linear BWM}
	In this section, we first derive an analytical expression for the optimal weights and examine the model's sensitivity to input data variation. We then compute CI to estimate the consistency of pairwise comparisons. Finally, we present numerical examples to demonstrate and validate the proposed framework.
	\subsection{Calculation of Weights}
	In this subsection, we derive a closed-form solution for the optimal weights and conduct sensitivity analysis of the model to data variations.\\\\
	\textbf{Notations:} 
	\begin{itemize}
		\item Let $\mathcal{W}=\left\{W=\{w_1,w_2,\ldots,w_n\}: \displaystyle\sum_{i=1}^{n}w_i=1 \text{ and } w_i\geq 0 \text{ for all } i\in C\right\}$.
		\item For $W\in\mathcal{W}$, let 
		\begin{equation}\label{4deviation_weight_set}
			\epsilon_W=\max\{|w_b-a_{bi}\times w_i|, |w_i-a_{iw}\times w_w|,|w_b-a_{bw}\times w_w|: i\in D\}.
		\end{equation}
		\item For $\lambda\geq 0$, let 
		\begin{equation}\label{4partition}
			\mathcal{W}_\lambda=\{W\in\mathcal{W}: \epsilon_W=\lambda\times w_w\}.
		\end{equation}
		Clearly, $\mathcal{W}=\displaystyle\bigcup_{\lambda\geq 0}\mathcal{W}_\lambda.$
		\item For $A=(A_b,A_w)$, let
		\begin{eqnarray*}
			&&D_1=\{i\in D: a_{bi}\times a_{iw}<a_{bw}\},\\
			&&D_2=\{i\in D: a_{bi}\times a_{iw}>a_{bw}\},\\
			&&D_3=\{i\in D: a_{bi}\times a_{iw}=a_{bw}\},
		\end{eqnarray*}
		\begin{equation}\label{4CV}
			\epsilon_i= \frac{|a_{bi}\times a_{iw} - a_{bw}|}{a_{bi}+2}, \quad 
			\epsilon_{i,j}= \frac{|a_{bi}\times a_{iw} - a_{bj}\times a_{jw}|}{a_{bi}+a_{bj}+2},\quad i,j\in D \quad \text{and}
		\end{equation}
		\begin{equation}\label{4max_1}
			\eta_A=\max\{\epsilon_i,\epsilon_{i,j}: i,j\in D\}.
		\end{equation}		
		Note that equation \eqref{4CV} can be rewritten as
		\begin{equation}\label{4CV1}
			\epsilon_i= \begin{cases}
				\frac{a_{bw}-a_{bi}\times a_{iw}}{a_{bi}+2} \quad \text{if } i\in D_1,\\
				\frac{a_{bi}\times a_{iw} - a_{bw}}{a_{bi}+2} \quad \text{if } i\in D_2\cup D_3,				
			\end{cases}
		\end{equation}
		\begin{equation}\label{4CV2}
			\epsilon_{i,j}= \begin{cases}
				\frac{a_{bj}\times a_{jw} - a_{bi}\times a_{iw}}{a_{bi}+a_{bj}+2} \quad \text{if } a_{bi}\times a_{iw}<a_{bj}\times a_{jw},\\
				\frac{a_{bi}\times a_{iw} - a_{bj}\times a_{jw}}{a_{bi}+a_{bj}+2} \quad \text{if } a_{bi}\times a_{iw}\geq a_{bj}\times a_{jw}.				
			\end{cases}
		\end{equation}\\
		Also, it is easy to verify that if $i,j\in D_1\cup D_3$ or $i,j\in D_2\cup D_3$, then $\epsilon_{i,j}<\max\{\epsilon_i,\epsilon_j\}$. This gives
		\begin{equation}\label{4max}
			\eta_A=\max\{\epsilon_i,\epsilon_{j,k}:i\in D, j\in D_1, k\in D_2\}.
		\end{equation} 
Whenever there is no ambiguity regarding $A$, we denote $\eta_A$ simply as $\eta$.
\end{itemize}
\begin{proposition}\label{4lower1}
	Let $W\in\mathcal{W}$, let $\epsilon_i$ and $\epsilon_{i,j}$ be as in equation \eqref{4CV}, and let $\epsilon_{W}$ be as in equation \eqref{4deviation_weight_set}. Then $\epsilon_i\times w_w,\epsilon_{i,j}\times w_w\leq\epsilon$ for all $i,j\in D$.
\end{proposition}
\begin{proof}
	Let $w_b-a_{bi}\times w_i= \eta_{bi}$, $w_i-a_{iw}\times w_w= \eta_{iw}$ and $w_b-a_{bw}\times w_w=\eta_{bw}$ for all $i\in D$. Thus, $|\eta_{bi}|,|\eta_{iw}|,|\eta_{bw}|\leq \epsilon$. Now, $w_b= a_{bi}\times w_i+\eta_{bi}$ and $w_i= a_{iw}\times w_w+\eta_{iw}$ gives $w_b= a_{bi}\times a_{iw}\times w_w+a_{bi}\times\eta_{iw}+\eta_{bi}$.
	Since we have $w_b= a_{bw}\times w_w+\eta_{bw}$, we get
	\begin{equation}\label{4deviation}
		a_{bw}\times w_w+\eta_{bw}=a_{bi}\times a_{iw}\times w_w+a_{bi}\times\eta_{iw}+\eta_{bi},
	\end{equation}
	where $i\in D$.\\\\
	Fix $i\in D$. If $i\in D_3$, then by equation \eqref{4CV1}, $\epsilon_i=0$. So, $\epsilon_i\times w_w= 0\leq \epsilon$. Now, assume that $i\in D_1$. Suppose, if possible, $|\eta_{bi}|,|\eta_{iw}|,|\eta_{bw}|<\epsilon_i \times w_w$, i.e., $\eta_{bi},\eta_{iw},\eta_{bw}\in (-\epsilon_i \times w_w,\epsilon_i \times w_w)$. Then, by equation \eqref{4deviation}, we get $a_{bw}\times w_w-\epsilon_i \times w_w<a_{bi}\times a_{iw}\times w_w+a_{bi}\times\epsilon_i\times w_w+\epsilon_i\times w_w$, which is contradiction as by equation \eqref{4CV1}, we have $a_{bw}\times w_w-\epsilon_i \times w_w=a_{bi}\times a_{iw}\times w_w+a_{bi}\times\epsilon_i\times w_w+\epsilon_i\times w_w$. So, $\epsilon_i \times w_w$ is less than or equal to at least one of $|\eta_{bi}|$, $|\eta_{iw}|$ and $|\eta_{bw}|$. Consequently, we have $\epsilon_i \times w_w\leq\epsilon$. Similar argument can be given if $i\in D_2$.\\\\
	Fix $i,j\in D$. If $a_{bi}\times a_{iw}=a_{bj}\times a_{jw}$, then by equation \eqref{4CV2}, $\epsilon_{i,j}=0$. So, $\epsilon_{i,j}\times w_W= 0\leq \epsilon$. If $a_{bi}\times a_{iw}\neq a_{bj}\times a_{jw}$, then without loss of generality, we may assume that $a_{bi}\times a_{iw}<a_{bj}\times a_{jw}$. Suppose, if possible, $|\eta_{bi}|,|\eta_{iw}|,|\eta_{bj}|,|\eta_{jw}|<\epsilon_{i,j} \times w_w$, i.e., $\eta_{bi},\eta_{iw},\eta_{bj},\eta_{jw}\in (-\epsilon_{i,j} \times w_w,\epsilon_{i,j} \times w_w)$. From equation \eqref{4deviation}, we have $a_{bi}\times a_{iw}\times w_w+a_{bi}\times\eta_{iw}+\eta_{bi}=a_{bj}\times a_{jw}\times w_w+a_{bj}\times\eta_{jw}+\eta_{bj}$. This gives $a_{bi}\times a_{iw}\times w_w+a_{bi}\times\epsilon_{i,j}\times w_w+\epsilon_{i,j}\times w_w>a_{bj}\times a_{jw}\times w_w-a_{bj}\times\epsilon_{i,j}\times w_w-\epsilon_{i,j}\times w_w$, which is contradiction as by equation \eqref{4CV2}, we have $a_{bi}\times a_{iw}\times w_w+a_{bi}\times\epsilon_{i,j}\times w_w+\epsilon_{i,j}\times w_w=a_{bj}\times a_{jw}\times w_w-a_{bj}\times\epsilon_{i,j}\times w_w-\epsilon_{i,j}\times w_w$. So, $\epsilon_{i,j} \times w_w$ is less than or equal to at least one of $|\eta_{bi}|$, $|\eta_{iw}|$, $|\eta_{bj}|$ and $|\eta_{jw}|$. Consequently, we have $\epsilon_{i,j} \times w_w\leq\epsilon$. Hence the proof. 
\end{proof}
Proposition \ref{4lower1} implies that $\mathcal{W}_\lambda=\phi$ for $\lambda<\eta$. So, we get $\mathcal{W}=\displaystyle\bigcup_{\lambda\geq \eta}\mathcal{W}_\lambda$.
\begin{proposition}\label{4non_empty}
	Let $\eta$ and $\mathcal{W}_\eta$ be as in equation \eqref{4max} and \eqref{4partition} respectively. Then $\mathcal{W}_\eta \neq \phi$.
\end{proposition}
\begin{proof}
	First assume that $\eta=\epsilon_{i_0}$ for some $i_0\in D_1$. Consider $(\tilde{A}_b,\tilde{A}_w)$ defined as
\begin{equation}\label{4modification}
	\begin{split}
		\tilde{a}_{bw}= a_{bw}-\epsilon_{i_0},\ \tilde{a}_{bi}= a_{bi}+\frac{\frac{a_{bw}-\epsilon_{i_0}-a_{bi}\times a_{iw}}{a_{bi}+1}}{a_{iw}+\frac{a_{bw}-\epsilon_{i_0}-a_{bi}\times a_{iw}}{a_{bi}+1}}, \ \tilde{a}_{iw}= a_{iw}+\frac{a_{bw}-\epsilon_{i_0}-a_{bi}\times a_{iw}}{a_{bi}+1},\ i\in D.
	\end{split}
\end{equation}
	It is easy to check that $\tilde{a}_{bi},\tilde{a}_{iw},\tilde{a}_{bw}>0$ for all $i\in D$, and thus, $(\tilde{A}_b,\tilde{A}_w)$ is well-defined. Note that $(\tilde{A}_b,\tilde{A}_w)$ is consistent. So, by Theorem \ref{4accurate}, there exists $W\in \mathcal{W}$ satisfying
	\begin{eqnarray*}
		&&\frac{w_b}{w_i}=\tilde{a}_{bi},\quad \frac{w_i}{w_w}=\tilde{a}_{iw},\quad \frac{w_b}{w_w}=\tilde{a}_{bw}, \quad i\in D,\\
		&&w_1+w_2+\ldots+w_n=1.
	\end{eqnarray*}
	This gives $w_b=\tilde{a}_{bi}\times w_i = \bigg(a_{bi}+\frac{\frac{a_{bw}-\epsilon_{i_0}-a_{bi}\times a_{iw}}{a_{bi}+1}}{a_{iw}+\frac{a_{bw}-\epsilon_{i_0}-a_{bi}\times a_{iw}}{a_{bi}+1}}\bigg)\times w_i = a_{bi}\times w_i + \bigg(\frac{\frac{a_{bw}-\epsilon_{i_0}-a_{bi}\times a_{iw}}{a_{bi}+1}}{a_{iw}+\frac{a_{bw}-\epsilon_{i_0}-a_{bi}\times a_{iw}}{a_{bi}+1}}\bigg) \times w_i$. Now, $w_i=\tilde{a}_{iw} \times w_w= \bigg(a_{iw}+\frac{a_{bw}-\epsilon_{i_0}-a_{bi}\times a_{iw}}{a_{bi}+1}\bigg) \times w_w$ implies $|w_b-a_{bi}\times w_i|=\bigg|\frac{a_{bw}-\epsilon_{i_0}-a_{bi}\times a_{iw}}{a_{bi}+1}\bigg|\times w_w$. Similarly, we get $|w_i-a_{iw}\times w_w|=\bigg|\frac{a_{bw}-\epsilon_{i_0}-a_{bi}\times a_{iw}}{a_{bi}+1}\bigg|\times w_w$ and $|w_b-a_{bw}\times w_w|=\epsilon_{i_0}\times w_w$.\\\\
	Now, to prove $\mathcal{W}_\eta(=\mathcal{W}_{\epsilon_{i_0}})\neq \phi$, it is sufficient to prove $\bigg|\frac{a_{bw}-\epsilon_{i_0}-a_{bi}\times a_{iw}}{a_{bi}+1}\bigg|\leq \epsilon_{i_0}$ for all $i\in D$. Fix $i\in D$. Then there are two possibilities.
	\begin{enumerate}[(i)]
		\item $a_{bi}\times a_{iw}< a_{bw}-\epsilon_{i_0}$\\\\
		In this case, we have $i\in D_1$. So, by equation \eqref{4CV1}, $\epsilon_i=\frac{a_{bi}\times a_{iw} - a_{bw}}{a_{bi}+2}$. This implies $a_{bi}\times a_{iw}+(a_{bi}+1)\times\epsilon_i=a_{bw}-\epsilon_i$. We also have $\bigg|\frac{a_{bw}-\epsilon_{i_0}-a_{bi}\times a_{iw}}{a_{bi}+1}\bigg|=\frac{a_{bw}-\epsilon_{i_0}-a_{bi}\times a_{iw}}{a_{bi}+1}$ as $\frac{a_{bw}-\epsilon_{i_0}-a_{bi}\times a_{iw}}{a_{bi}+1}> 0$. Since $(\tilde{A}_b,\tilde{A}_w)$ is consistent, $\tilde{a}_{bi}\times \tilde{a}_{iw}=\tilde{a}_{bw}$ holds. So, by equation \eqref{4modification}, we get $$\bigg(a_{bi}+\frac{\frac{a_{bw}-\epsilon_{i_0}-a_{bi}\times a_{iw}}{a_{bi}+1}}{a_{iw}+\frac{a_{bw}-\epsilon_{i_0}-a_{bi}\times a_{iw}}{a_{bi}+1}}\bigg)\times \bigg(a_{iw}+\frac{a_{bw}-\epsilon_{i_0}-a_{bi}\times a_{iw}}{a_{bi}+1}\bigg)=a_{bw}-\epsilon_{i_0}.$$ This gives $a_{bi}\times a_{iw}+(a_{bi}+1)\times \bigg(\frac{a_{bw}-\epsilon_{i_0}-a_{bi}\times a_{iw}}{a_{bi}+1}\bigg)=a_{bw}-\epsilon_{i_0}$. Now, $\epsilon_i\leq \epsilon_{i_0}$ implies $a_{bw}-\epsilon_{i_0}\leq a_{bw}-\epsilon_i$, i.e., $$a_{bi}\times a_{iw}+(a_{bi}+1)\times \bigg(\frac{a_{bw}-\epsilon_{i_0}-a_{bi}\times a_{iw}}{a_{bi}+1}\bigg)\leq a_{bi}\times a_{iw}+(a_{bi}+1)\times\epsilon_i.$$ Thus, $\frac{a_{bw}-\epsilon_{i_0}-a_{bi}\times a_{iw}}{a_{bi}+1}\leq \epsilon_{i_0}$.\\
		\item $a_{bi}\times a_{iw}\geq a_{bw}-\epsilon_{i_0}$\\\\
		From equation \eqref{4CV1}, we have $a_{bi_0}\times a_{i_0w}+(a_{bi_0}+1)\times\epsilon_{i_0}=a_{bw}-\epsilon_{i_0}$. So, in this case, $a_{bi}\times a_{iw}>a_{bi_0}\times a_{i_0w}$. Thus, by equation \eqref{4CV2}, $a_{bi}\times a_{iw}-(a_{bi}+1)\times\epsilon_{i,i_0}=a_{bi_0}\times a_{i_0w}+(a_{bi_0}+1)\times \epsilon_{i,i_0}$. Also, $\bigg|\frac{a_{bw}-\epsilon_{i_0}-a_{bi}\times a_{iw}}{a_{bi}+1}\bigg|=-\bigg(\frac{a_{bw}-\epsilon_{i_0}-a_{bi}\times a_{iw}}{a_{bi}+1}\bigg)$ as $\frac{a_{bw}-\epsilon_{i_0}-a_{bi}\times a_{iw}}{a_{bi}+1}\leq 0$. As discussed in possibility (i), $a_{bi}\times a_{iw}+(a_{bi}+1)\times \bigg(\frac{a_{bw}-\epsilon_{i_0}-a_{bi}\times a_{iw}}{a_{bi}+1}\bigg)=a_{bw}-\epsilon_{i_0}$ holds. So, we get $$a_{bi}\times a_{iw}+(a_{bi}+1)\times \bigg(\frac{a_{bw}-\epsilon_{i_0}-a_{bi}\times a_{iw}}{a_{bi}+1}\bigg)=a_{bi_0}\times a_{i_0w}+(a_{bi_0}+1)\times\epsilon_{i_0}.$$ Now, $\epsilon_{i,i_0}\leq \epsilon_{i_0}$ implies $a_{bi_0}\times a_{i_0w}+(a_{bi_0}+1)\times \epsilon_{i,i_0}\leq a_{bi_0}\times a_{i_0w}+(a_{bi_0}+1)\times\epsilon_{i_0}$, i.e., $$a_{bi}\times a_{iw}-(a_{bi}+1)\times\epsilon_{i,i_0}\leq a_{bi}\times a_{iw}+(a_{bi}+1)\times \bigg(\frac{a_{bw}-\epsilon_{i_0}-a_{bi}\times a_{iw}}{a_{bi}+1}\bigg).$$ This gives $0<-\bigg(\frac{a_{bw}-\epsilon_{i_0}-a_{bi}\times a_{iw}}{a_{bi}+1}\bigg)\leq \epsilon_{i,i_0}\leq \epsilon_{i_0}$.
	\end{enumerate}
	Similar argument can be given if $\eta=\epsilon_{j_0}$ for some $j_0\in D_2$, or $\eta=\epsilon_{i_0,j_0}$ for some $i_0\in D_1$ and $j_0\in D_2$. Hence the proof.
\end{proof}
\begin{proposition}\label{4fix}
	Let $W\in\mathcal{W}_\eta$. Then 
	\begin{equation}\label{4fix_values}
		\begin{split}
			\begin{cases}
				\begin{cases}
					w_b&=\bigg(a_{bi_0}+\frac{\epsilon_{i_0}}{a_{i_0w}+\epsilon_{i_0}}\bigg)\times w_{i_0}\\
					&=(a_{bw}-\epsilon_{i_0})\times w_w,\\
					w_{i_0}&=(a_{i_0w}+\epsilon_{i_0})\times w_w
				\end{cases}\quad\quad\quad\quad\quad\quad\text{if } \eta=\epsilon_{i_0} \text{ for some }i_0\in D_1,\\
				\begin{cases}
					w_b&=\bigg(a_{bj_0}-\frac{\epsilon_{j_0}}{a_{j_0w}-\epsilon_{j_0}}\bigg)\times w_{j_0}\\
					&=(a_{bw}+\epsilon_{j_0})\times w_w,\\
					w_{j_0}&=(a_{j_0w}-\epsilon_{j_0})\times w_w 
				\end{cases}\quad\quad\quad\quad\quad\quad\text{if } \eta=\epsilon_{j_0} \text{ for some }j_0\in D_2,\\
				\begin{cases}
					w_b&=\bigg(a_{bi_0}+\frac{\epsilon_{i_0,j_0}}{a_{i_0w}+\epsilon_{i_0,j_0}}\bigg)\times w_{i_0}\\ &=\bigg(a_{bj_0}-\frac{\epsilon_{i_0,j_0}}{a_{j_0w}-\epsilon_{i_0,j_0}}\bigg)\times w_{j_0}\\
					&=(a_{bi_0}\times a_{i_0w}+(a_{bi_0}+1)\times \epsilon_{i_0,j_0})\times w_w\\
					&=(a_{bj_0}\times a_{j_0w}-(a_{bj_0}+1)\times \epsilon_{i_0,j_0})\times w_w,\\
					w_{i_0}&=(a_{i_0w}+\epsilon_{i_0,j_0})\times w_w,\\
					w_{j_0}&=(a_{j_0w}-\epsilon_{i_0,j_0})\times w_w\\
				\end{cases}\text{if } \eta=\epsilon_{i_0,j_0} \text{ for some }i_0\in D_1 \text{ and } j_0\in D_2.
			\end{cases}
		\end{split}
	\end{equation}
\end{proposition}
\begin{proof}
	Since all three cases are similar, we prove the result for the first case only. Let $\eta=\epsilon_{i_0}$ for some $i_0\in D_1$. As discussed in Proposition \ref{4lower1}, we get $\eta_{bi_0},\eta_{i_0w},\eta_{bw} \in [-\epsilon_{i_0}\times w_w,\epsilon_{i_0}\times w_w]$ such that $w_b=a_{bi_0}\times w_{i_0}+\eta_{bi_0}=a_{bw}\times w_w+\eta_{bw}$ and $w_{i_0}=a_{i_0w}\times w_w+\eta_{bw}$, and thus, $a_{bw}\times w_w+\eta_{bw}=a_{bi_0}\times a_{i_0w}\times w_w+a_{bi_0}\times\eta_{i_0w}+\eta_{bi_0}$. Also, by equation \eqref{4CV1}, we have $a_{bw}\times w_w-\epsilon_{i_0}\times w_w=a_{bi_0}\times a_{i_0w}\times w_w+a_{bi_0}\times \epsilon_{i_0}\times w_w+\epsilon_{i_0}\times w_w$. So, if either $\eta_{bi_0}$ or $\eta_{i_0w}$ is strictly less than $\epsilon_{i_0}\times w_w$, then we get $a_{bi_0}\times a_{i_0w}\times w_w+a_{bi_0}\times\eta_{i_0w}+\eta_{bi_0}<a_{bi_0}\times a_{i_0w}\times w_w+a_{bi_0}\times \epsilon_{i_0}\times w_w+\epsilon_{i_0}\times w_w$, which implies $a_{bw}\times w_w+\eta_{bw}<a_{bw}\times w_w-\epsilon_{i_0}\times w_w$, i.e, $\eta_{bw}<-\epsilon_{i_0}\times w_w$, which is contradiction as $\eta_{bw} \in [-\epsilon_{i_0}\times w_w,\epsilon_{i_0}\times w_w]$. Therefore, $\eta_{bi_0}=\eta_{i_0w}=\epsilon_{i_0}\times w_w$, and thus, $\eta_{bw}=-\epsilon_{i_0}\times w_w$. Now, $w_b=a_{bi_0}\times w_{i_0}+\eta_{bi_0}=a_{bw}\times w_w+\eta_{bw}$ and $w_{i_0}=a_{i_0w}\times w_w+\eta_{bw}$ gives $w_b=\bigg(a_{bi_0}+\frac{\epsilon_{i_0}}{a_{i_0w}+\epsilon_{i_0}}\bigg)\times w_{i_0}=(a_{bw}-\epsilon_{i_0})\times w_w$ and $w_{i_0}=(a_{i_0w}+\epsilon_{i_0})\times w_w$. Hence the proof.
\end{proof}
\begin{lemma}\label{4inequality}
	Let $\epsilon_i$ and $\epsilon_{i,j}$ be as in equation \eqref{4CV}, and let $\eta$ be as in equation \eqref{4max}. Then the following statements hold.
	\begin{enumerate}
		\item If $\eta=\epsilon_{i_0}$ for some $i_0 \in D_1$, then $\max\biggl\{a_{iw}-\epsilon_{i_0},\frac{a_{bw}-2\epsilon_{i_0}}{a_{bi}}\biggr\}\leq \min\biggl\{a_{iw}+\epsilon_{i_0},\frac{a_{bw}}{a_{bi}}\biggr\}$ for all $i\in D$.
		\item If $\eta=\epsilon_{j_0}$ for some $j_0 \in D_2$, then $\max\biggl\{a_{iw}-\epsilon_{j_0},\frac{a_{bw}}{a_{bi}}\biggr\}\leq \min\biggl\{a_{iw}+\epsilon_{j_0},\frac{a_{bw}+2\epsilon_{j_0}}{a_{bi}}\biggr\}$ for all $i\in D$.
		\item If $\eta=\epsilon_{i_0,j_0}$ for some $i_0 \in D_1$ and $j_0\in D_2$, then $\max\biggl\{a_{iw}-\epsilon_{i_0,j_0},\frac{a_{bi_0}\times a_{i_0w}+a_{bi_0}\times \epsilon_{i_0,j_0}}{a_{bi}}\biggr\}\leq \min\biggl\{a_{iw}+\epsilon_{i_0,j_0},\frac{a_{bi_0}\times a_{i_0w}+(a_{bi_0}+2)\times \epsilon_{i_0,j_0}}{a_{bi}}\biggr\}$ for all $i\in D$.
	\end{enumerate}
\end{lemma}
\begin{proof}
	 Here, we prove only the first statement, as the proofs of all three statements are similar. Let $W\in\mathcal{W}_{\epsilon_{i_0}}$. So, we have $|w_b-a_{bi}\times w_i|\leq \epsilon_{i_0}\times w_w$, i.e., $-\epsilon_{i_0}\times w_w \leq w_b-a_{bi}\times w_i\leq \epsilon_{i_0}\times w_w$ for all $i\in D$. Now, Proposition \ref{4fix} gives $w_b=(a_{bw}-\epsilon_{i_0})\times w_w$. Thus, we get $$\biggl(\frac{a_{bw}-2\epsilon_{i_0}}{a_{bi}}\biggr)\times w_w\leq w_i\leq \biggl(\frac{a_{bw}}{a_{bi}}\biggr)\times w_w.$$ Similarly, $|w_i-a_{iw}\times w_w|\leq \epsilon_{i_0}\times w_w$ implies $$(a_{iw}-\epsilon_{i_0})\times w_w\leq w_i\leq (a_{iw}+\epsilon_{i_0})\times w_w.$$ Combining these inequalities, we get 
	 \begin{equation}\label{4inequality1}
	 	\max\biggl\{a_{iw}-\epsilon_{i_0},\frac{a_{bw}-2\epsilon_{i_0}}{a_{bi}}\biggr\}\times w_w\leq w_i \leq \min\biggl\{a_{iw}+\epsilon_{i_0},\frac{a_{bw}}{a_{bi}}\biggr\}\times w_w.
	 \end{equation}
	 Thus, $\max\biggl\{a_{iw}-\epsilon_{i_0},\frac{a_{bw}-2\epsilon_{i_0}}{a_{bi}}\biggr\}\leq \min\biggl\{a_{iw}+\epsilon_{i_0},\frac{a_{bw}}{a_{bi}}\biggr\}$. Hence the proof.
\end{proof}
\begin{theorem}\label{4exact_obj}
	let $\epsilon_i$ and $\epsilon_{i,j}$ be as in equation \eqref{4CV}, let $\epsilon_{W}$ and $\eta$ be as in equation \eqref{4deviation_weight_set} and \eqref{4max} respectively, and let $\epsilon^*$ be the optimal objective value of problem \eqref{4optimization_2}. Then the following statements hold.
	\begin{enumerate}
		\item If $\eta=\epsilon_{i_0}$ for some $i_0 \in D_1$, then $W=\{w_1,w_2,\ldots,w_n\}$ defined as
		\begin{equation}\label{4optimal1}
			\begin{split}
				&w_b=\frac{a_{bw}-\epsilon_{i_0}}{\sigma},\quad 
				w_{i_0}=\frac{a_{i_0w}+\epsilon_{i_0}}{\sigma}, \quad 
				w_w=\frac{1}{\sigma},\\
				&w_i=\frac{\min\biggl\{a_{iw}+\epsilon_{i_0},\frac{a_{bw}}{a_{bi}}\biggr\}}{\sigma}\text{ for } i\in D\setminus\{i_0\},
			\end{split}
		\end{equation}
		where $\sigma=1+a_{bw}+a_{i_0w}+\displaystyle\sum_{\substack{i\in D\\i\neq i_0}}\min\biggl\{a_{iw}+\epsilon_{i_0},\frac{a_{bw}}{a_{bi}}\biggr\}$ is the optimal weight set, and consequently, $\epsilon^*=\epsilon_W=\frac{\epsilon_{i_0}}{\sigma}$.
		\item If $\eta=\epsilon_{j_0}$ for some $j_0\in D_2$, then $W=\{w_1,w_2,\ldots,w_n\}$ defined as
		\begin{equation}\label{4optimal2}
			\begin{split}
				&w_b=\frac{a_{bw}+\epsilon_{j_0}}{\sigma},\quad 
				w_{j_0}=\frac{a_{j_0w}-\epsilon_{j_0}}{\sigma}, \quad
				w_w=\frac{1}{\sigma},\\
				&w_i=\frac{\min\biggl\{a_{iw}+\epsilon_{j_0},\frac{a_{bw}+2\epsilon_{j_0}}{a_{bi}}\biggr\}}{\sigma}\text{ for } i\in D\setminus\{j_0\},
			\end{split}
		\end{equation}
		where $\sigma=1+a_{bw}+a_{j_0w}+\displaystyle\sum_{\substack{i\in D\\i\neq j_0}}\min\biggl\{a_{iw}+\epsilon_{j_0},\frac{a_{bw}+2\epsilon_{j_0}}{a_{bi}}\biggr\}$ is the optimal weight set, and consequently, $\epsilon^*=\epsilon_W=\frac{\epsilon_{j_0}}{\sigma}$.
		\item If $\eta=\epsilon_{i_0,j_0}$ for some $i_0\in D_1$ and $j_0\in D_2$, then $W=\{w_1,w_2,\ldots,w_n\}$ defined as
		\begin{equation}\label{4optimal3}
			\begin{split}
				&w_b=\frac{a_{bi_0}\times a_{i_0w}+(a_{bi_0}+1)\times \epsilon_{i_0,j_0}}{\sigma},\quad
				w_{i_0}=\frac{a_{i_0w}+\epsilon_{i_0,j_0}}{\sigma},\quad
				w_{j_0}=\frac{a_{j_0w}-\epsilon_{i_0,j_0}}{\sigma},\\ 	
				&w_w=\frac{1}{\sigma},\quad
				w_i=\frac{\min\biggl\{a_{iw}+\epsilon_{i_0,j_0},\frac{a_{bi_0}\times a_{i_0w}+(a_{bi_0}+2)\times \epsilon_{i_0,j_0}}{a_{bi}}\biggr\}}{\sigma}\text{ for } i\in D\setminus\{i_0,j_0\},
			\end{split}
		\end{equation} 
		$\begin{aligned}
			\text{where }\sigma=1+a_{i_0w}+a_{j_0w}+a_{bi_0}\times a_{i_0w}+(a_{bi_0}+1)\times \epsilon_{i_0,j_0}\quad\quad\quad\quad\quad\quad\quad\quad\quad\quad\quad\quad\quad&&\\+\displaystyle\sum_{\substack{i\in D\\i\neq i_0,j_0}}\min\biggl\{a_{iw}+\epsilon_{i_0,j_0},\frac{a_{bi_0}\times a_{i_0w}+(a_{bi_0}+2)\times \epsilon_{i_0,j_0}}{a_{bi}}\biggr\}&&
		\end{aligned}$\\
		is the optimal weight set, and consequently, $\epsilon^*=\epsilon_W=\frac{\epsilon_{i_0,j_0}}{\sigma}$.
		\end{enumerate}
\end{theorem}
\begin{proof}
	Here, we present the proof of the first statement only, as the proofs of the other two statements follow a similar reasoning. \\\\
	Here, $$|w_b-a_{bi_0}\times w_{i_0}|=\bigg|\frac{a_{bw}-\epsilon_{i_0}-a_{bi_0}\times (a_{i_0w}+\epsilon_{i_0})}{\sigma}\bigg|=\bigg|\frac{a_{bw}-a_{bi_0}\times a_{i_0w}-(a_{bi_0}+1)\times\epsilon_{i_0}}{\sigma}\bigg|.$$ Using equation \eqref{4CV1}, we get $|w_b-a_{bi_0}\times w_{i_0}|=\big|\frac{\epsilon_{i_0}}{\sigma}\big|=\frac{\epsilon_{i_0}}{\sigma}$. Also, $$|w_{i_0}-a_{i_0w}\times w_w|=\bigg|\frac{a_{i_0w}+\epsilon_{i_0}-a_{i_0w}\times 1}{\sigma}\bigg|=\bigg|\frac{\epsilon_{i_0}}{\sigma}\bigg|=\frac{\epsilon_{i_0}}{\sigma} \text{ and }$$  $$|w_{b}-a_{bw}\times w_w|=\bigg|\frac{a_{bw}-\epsilon_{i_0}-a_{bw}\times 1}{\sigma}\bigg|=\bigg|-\frac{\epsilon_{i_0}}{\sigma}\bigg|=\frac{\epsilon_{i_0}}{\sigma}.$$ Now, fix $i\in D\setminus\{i_0\}$. Then $$w_b-a_{bi}\times w_i=\frac{a_{bw}-\epsilon_{i_0}-a_{bi}\times \min\bigl\{a_{iw}+\epsilon_{i_0},\frac{a_{bw}}{a_{bi}}\bigr\}}{\sigma}.$$ This, along with statement \textit{1} of Lemma \ref{4inequality}, implies that
	\begin{eqnarray*}
		-\frac{\epsilon_{i_0}}{\sigma}=\frac{a_{bw}-\epsilon_{i_0}-a_{bi}\times \bigl(\frac{a_{bw}}{a_{bi}}\bigr)}{\sigma}\leq w_b-a_{bi}\times w_i&\leq& \frac{a_{bw}-\epsilon_{i_0}-a_{bi}\times \max\bigl\{\frac{a_{bw}-2\epsilon_{i_0}}{a_{bi}},a_{iw}-\epsilon_{i_0}\bigr\}}{\sigma}\\
		&\leq& \frac{a_{bw}-\epsilon_{i_0}-a_{bi}\times \bigl(\frac{a_{bw}-2\epsilon_{i_0}}{a_{bi}}\bigr)}{\sigma}=\frac{\epsilon_{i_0}}{\sigma}.
	\end{eqnarray*}
	Thus, $|w_b-a_{bi}\times w_i|\leq \frac{\epsilon_{i_0}}{\sigma}$. Similarly, we get $|w_i-a_{iw}\times w_w|\leq \frac{\epsilon_{i_0}}{\sigma}$. Therefore, $\epsilon_W=\frac{\epsilon_{i_0}}{\sigma}$. Now, it remains to show that $W$ is the optimal weight set.\\\\
	Let $W'=\{w_1',w_2',\ldots,w_n'\}\in \mathcal{W}$. So, $W'\in \mathcal{W}_\lambda$ for some $\lambda\geq \epsilon_{i_0}$. Then there are two possibilities.
	\begin{enumerate}[(i)]
		\item $\lambda=\epsilon_{i_0}$\\\\
		In this case, by Proposition \ref{4fix}, we have $w_b'=(a_{bw}-\epsilon_{i_0})\times w_w'$ and $w_{i_0}'=(a_{i_0w}+\epsilon_{i_0})\times w_w'$, and from equation \eqref{4inequality1}, we have $w_i'\leq \min\bigl\{\frac{a_{bw}}{a_{bi}},a_{iw}+\epsilon_{i_0}\bigr\}\times w_w'$ for all $i\in D\setminus\{i_0\}$. Now, $\displaystyle\sum_{i=1}^{n}w_i'=1$ implies that $$\bigg(1+a_{bw}+a_{i_0w}+\displaystyle\sum_{\substack{i\in D\\i\neq i_0}}\min\biggl\{a_{iw}+\epsilon_{i_0},\frac{a_{bw}}{a_{bi}}\biggr\}\bigg)\times w_w'\geq 1,$$ i.e., $w_w'\geq \frac{1}{\sigma}=w_w$. Now, $W'\in \mathcal{W}_{\epsilon_{i_0}}$ implies that $\epsilon_{W'}=\epsilon_{i_0}\times w_w'\geq \epsilon_{i_0}\times w_w=\epsilon_W$.\\\\
		Observe that if $\epsilon_{W'}=\epsilon_{W}$, then $w_w'=w_w=\frac{1}{\sigma}$, which is true only if $$w_i'= \min\biggl\{\frac{a_{bw}}{a_{bi}},a_{iw}+\epsilon_{i_0}\biggr\}\times w_w'=\frac{\min\biggl\{\frac{a_{bw}}{a_{bi}},a_{iw}+\epsilon_{i_0}\biggr\}}{\sigma}=w_i$$ for all $i\in D\setminus\{i_0\}$. Also, $$w_b'=(a_{bw}-\epsilon_{i_0})\times w_w'=\frac{a_{bw}-\epsilon_{i_0}}{\sigma}=w_b \quad \text{and}\quad w_{i_0}'=(a_{i_0w}+\epsilon_{i_0})\times w_w'=\frac{a_{i_0w}+\epsilon_{i_0}}{\sigma}.$$ So, if $\epsilon_{W'}=\epsilon_{W}$, then $W'=W$.\\
		\item $\lambda>\epsilon_{i_0}$\\\\
		Let $\delta=\lambda-\epsilon_{i_0}$. So, $\delta>0$. Let 
		\begin{equation}\label{4deviation1}
			w_b'-a_{bi}\times w_i'= \eta_{bi}',\quad w_i'-a_{iw}\times w_w'= \eta_{iw}',\quad w_b'-a_{bw}\times w_w'=\eta_{bw}',\quad \text{where } i\in D.
		\end{equation}
		Since $W'\in \mathcal{W}_\lambda$, we get $|\eta_{bi}'|,|\eta_{iw}'|,|\eta_{bw}'|\leq \lambda \times w_w'$. From equation \eqref{4deviation1}, it follows that $$a_{bw}\times w_w'+\eta_{bw}'=a_{bi_0}\times a_{i_0w}\times w_w'+a_{bi_0}\times\eta_{i_0w}'+\eta_{bi_0}'.$$ Now, $\eta_{bi_0}'\leq |\eta_{bi_0}'|\leq \lambda \times w_w'=(\epsilon_{i_0}+\delta)\times w_w'$ and $\eta_{i_0w}'\leq |\eta_{i_0w}'|\leq \lambda \times w_w'=(\epsilon_{i_0}+\delta)\times w_w'$ imply that $$a_{bw}\times w_w'+\eta_{bw}'\leq a_{bi_0}\times a_{i_0w}\times w_w'+a_{bi_0}\times(\epsilon_{i_0}+\delta)\times w_w'+(\epsilon_{i_0}+\delta)\times w_w'.$$ Using equation \eqref{4CV1}, we get $$a_{bw}\times w_w'+\eta_{bw}'\leq a_{bw}\times w_w'-\epsilon_{i_0}\times w_w'+(a_{bi_0}+1)\delta\times w_w',$$ i.e., $\eta_{bw}'\leq (a_{bi_0}+1-\epsilon_{i_0})\delta\times w_w'$. Now, equation \eqref{4deviation1} gives 
		\begin{eqnarray*}
			&&w_b'\leq (a_{bw}-\epsilon_{i_0}+(a_{bi_0}+1)\delta)\times w_w',\quad w_{i_0}'\leq (a_{i_0w}+\epsilon_{i_0}+\delta)\times w_w',\\
			&&w_i'\leq (a_{iw}+\epsilon_{i_0}+\delta)\times w_w',\quad w_i'\leq\frac{w_b'+\lambda\times w_w'}{a_{bi}}\leq \frac{a_{bw}+(a_{bi_0}+2)\delta}{a_{bi}}\times w_w'. 
		\end{eqnarray*}
		So, $$w_i'\leq \min\biggl\{a_{iw}+\epsilon_{i_0}+\delta,\frac{a_{bw}+(a_{bi_0}+2)\delta}{a_{bi}}\biggr\}\times w_w'$$ for all $i\in D\setminus \{i_0\}$. Now, $\displaystyle\sum_{i=1}^{n}w_i'=1$ implies that $$\bigg(1+a_{bw}+a_{i_0w}+(a_{bi_0}+2)\delta+\displaystyle\sum_{\substack{i\in D\\i\neq i_0}}\min\biggl\{a_{iw}+\epsilon_{i_0}+\delta,\frac{a_{bw}+(a_{bi_0}+2)\delta}{a_{bi}}\biggr\}\bigg)\times w_w'\geq 1.$$ Since $W'\in \mathcal{W}_{\lambda}$, we get $$\frac{\epsilon_{i_0}+\delta}{1+a_{bw}+a_{i_0w}+(a_{bi_0}+2)\delta+\displaystyle\sum_{\substack{i\in D\\i\neq i_0}}\min\biggl\{a_{iw}+\epsilon_{i_0}+\delta,\frac{a_{bw}+(a_{bi_0}+2)\delta}{a_{bi}}\biggr\}}\leq \lambda\times w_w'=\epsilon_{W'}.$$ This, along with the fact that $$f(x)=\frac{\epsilon_{i_0}+x}{1+a_{bw}+a_{i_0w}+(a_{bi_0}+2)x+\displaystyle\sum_{\substack{i\in D\\i\neq i_0}}\min\biggl\{a_{iw}+\epsilon_{i_0}+x,\frac{a_{bw}+(a_{bi_0}+2)x}{a_{bi}}\biggr\}},$$ where $x\in [0,\infty)$ is a strictly increasing function, gives $$\epsilon_W= \frac{\epsilon_{i_0}}{1+a_{bw}+a_{i_0w}+\displaystyle\sum_{\substack{i\in D\\i\neq i_0}}\min\biggl\{a_{iw}+\epsilon_{i_0},\frac{a_{bw}}{a_{bi}}\biggr\}}<\epsilon_{W'}.$$
	\end{enumerate}
	From the above discussion, it follows that $W$ defined by \eqref{4optimal1} is the optimal weight set. Hence the proof.
\end{proof}
\begin{corollary}\label{4sa}
	Let $A=(A_b,A_w)$ be a PCS, let $\epsilon_i$ and $\epsilon_{i,j}$ be as in equation \eqref{4CV}, and let $\eta_A$ be as in equation \eqref{4max}. Then the following statements hold.
	\begin{enumerate}
		\item If $\eta=\epsilon_{i_0}$ for some $i_0\in D_1$, then for any PCS $A'=(A_b',A_w')$ such that $a_{bi_0}'=a_{bi_0}$, $a_{i_0w}'=a_{i_0w}$, $a_{bw}'=a_{bw}$, $\eta_{A'}=\eta_A=\epsilon_{i_0}$ and $\min\biggl\{a_{iw}'+\epsilon_{i_0},\frac{a_{bw}'}{a_{bi}'}\biggr\}=\min\biggl\{a_{iw}+\epsilon_{i_0},\frac{a_{bw}}{a_{bi}}\biggr\}$ for all $i\in D\setminus\{i_0\}$, both the optimal weight set and the optimal objective value of problem \eqref{4optimization_2} coincide with those of $A$.
		\item If $\eta=\epsilon_{j_0}$ for some $j_0\in D_2$, then for any PCS $A'=(A_b',A_w')$ such that $a_{bj_0}'=a_{bj_0}$, $a_{j_0w}'=a_{j_0w}$, $a_{bw}'=a_{bw}$, $\eta_{A'}=\eta_A=\epsilon_{j_0}$ and $\min\biggl\{a_{iw}'+\epsilon_{j_0},\frac{a_{bw}'+2\epsilon_{j_0}}{a_{bi}'}\biggr\}=\min\biggl\{a_{iw}+\epsilon_{j_0},\frac{a_{bw}+2\epsilon_{j_0}}{a_{bi}}\biggr\}$ for all $i\in D\setminus\{j_0\}$, both the optimal weight set and the optimal objective value of problem \eqref{4optimization_2} coincide with those of $A$.
		\item If $\eta=\epsilon_{i_0,j_0}$ for some $i_0\in D_1$ and $j_0\in D_2$, then for any PCS $A'=(A_b',A_w')$ such that $a_{bi_0}'=a_{bi_0}$, $a_{i_0w}'=a_{i_0w}$, $a_{bj_0}'=a_{bj_0}$, $a_{j_0w}'=a_{j_0w}$, $a_{bw}'=a_{bw}$, $\eta_{A'}=\eta_A=\epsilon_{i_0,j_0}$ and $\min\biggl\{a_{iw}'+\epsilon_{i_0,j_0},\frac{a_{bi_0}'\times a_{i_0w}'+(a_{bi_0}'+2)\times \epsilon_{i_0,j_0}}{a_{bi}'}\biggr\}=\min\biggl\{a_{iw}+\epsilon_{i_0,j_0},\frac{a_{bi_0}\times a_{i_0w}+(a_{bi_0}+2)\times \epsilon_{i_0,j_0}}{a_{bi}}\biggr\}$ for all $i\in D\setminus\{i_0,j_0\}$, both the optimal weight set and the optimal objective value of problem \eqref{4optimization_2} coincide with those of $A$.
	\end{enumerate}
\end{corollary}
\hspace{-0.7cm}
Corollary \ref{4sa} reveals that in some instances, even after modifying the given comparison value, the optimal weight set remains unchanged, with the same corresponding $\epsilon^*$. This suggests that the linear BWM exhibits low sensitivity to data perturbations, which is one of its limitations.
\subsection{Consistency Analysis}
For any MCDM method, assessing consistency of the decision data is crucial as the outcome depends directly on these inputs. In BWM, this assessment is preformed using the Consistency Ratio (CR), defined as
\begin{equation}\label{4cr}
	\text{CR}=\frac{\epsilon^*}{\text{Consistency Index (CI)}},
\end{equation}
where CI $=\sup\left\{\epsilon^*: \epsilon^*\text{ is the optimal objective value of problem \eqref{4optimization_2} for some } (A_b,A_w) \text{ having} \right.\\ \left.\text{the given number of criteria } n \text{ and the given value of } a_{bw}\right\}$ \cite{rezaei2015best}. So, CI is a function of $n$ and $a_{bw}$. In this subsection, our objective is to compute CI$_{a_{bw}}(n)$ for the linear BWM.
\begin{theorem}\label{4ci}
	Let $(A_b,A_w)$ be a PCS with $n$ criteria and having $a_{bw}$ as the best-to-worst comparison value, let $\epsilon^*$ be the corresponding optimal objective value of problem \eqref{4optimization_2}, let $\epsilon_i$ and $\epsilon_{i,j}$ be as in equation \eqref{4CV}, and let $\eta$ be as in equation \eqref{4max}. Then the following statements hold.
	\begin{enumerate}
		\item For $n\geq 3$, if $\eta=\epsilon_{i_0}$ for some $i_0 \in D_1$, then $\epsilon^*\leq \frac{a_{bw}-1}{3(n-1+a_{bw})}$. Furthermore, $(A_b',A_w')$ defined as $$
		a_{1n}'=a_{bw},\ a_{12}'=a_{2n}'=1,\ a_{1i}'=a_{bw}, \ a_{in}'=1 \text{ for }i\in\{3,4,\ldots,n-1\}$$ is a PCS with $n$ criteria with $c_1$ as the best and $c_n$ as the worst criterion, having $a_{bw}$ as the best-to-worst comparison value, for which $\epsilon^*=\frac{a_{bw}-1}{3(n-1+a_{bw})}$.
		\item For $n\geq 3$, if $\eta=\epsilon_{j_0}$ for some $j_0 \in D_2$, then $\epsilon^*\leq \frac{a_{bw}(a_{bw}-1)}{2a_{bw}^2+(3n-4)a_{bw}+2}$. Furthermore, $(A_b',A_w')$ defined as $$
		a_{1n}'=a_{bw},\ a_{12}'=a_{2n}'=a_{bw},\ a_{1i}'=a_{bw}, \ a_{in}'=1 \text{ for }i\in\{3,4,\ldots,n-1\}$$ is a PCS with $n$ criteria with $c_1$ as the best and $c_n$ as the worst criterion, having $a_{bw}$ as the best-to-worst comparison value, for which $\epsilon^*=\frac{a_{bw}(a_{bw}-1)}{2a_{bw}^2+(3n-4)a_{bw}+2}$.
		\item For $n\geq 4$, if $\eta=\epsilon_{i_0,j_0}$ for some $i_0\in D_1$ and $j_0 \in D_2$, then\\ $\epsilon^*\leq \frac{a_{bw}^2-1}{3a_{bw}^2+6a_{bw}+7+(n-4)\min\{a_{bw}^2+a_{bw}+2,1+3a_{bw}\}}$. Furthermore, $(A_b',A_w')$ defined as $$
		a_{1n}'=a_{bw},\ a_{12}'=a_{2n}'=1,\ a_{13}'=a_{3n}'=a_{bw},\ a_{1i}'=a_{bw}, \ a_{in}'=1 \text{ for }i\in\{4,\ldots,n-1\}$$ is a PCS with $n$ criteria with $c_1$ as the best and $c_n$ as the worst criterion, having $a_{bw}$ as the best-to-worst comparison value, for which $\epsilon^*=\frac{a_{bw}^2-1}{3a_{bw}^2+6a_{bw}+7+(n-4)\min\{a_{bw}^2+a_{bw}+2,1+3a_{bw}\}}$.
	\end{enumerate}
\end{theorem}
\begin{proof}
	Here, we prove statement \textit{1} only, as statements \textit{2} and \textit{3} follow analogously.\\\\
	Since $\eta=\epsilon_{i_0}$ for $i_0 \in D_1$, by statement \textit{1} of Theorem \ref{4exact_obj}, we have $$\epsilon^*=\frac{\epsilon_{i_0}}{1+a_{bw}+a_{i_0w}+\displaystyle\sum_{\substack{i\in D\\i\neq i_0}}\min\biggl\{a_{iw}+\epsilon_{i_0},\frac{a_{bw}}{a_{bi}}\biggr\}}.$$
	For the Saaty scale, we have $a_{iw}\geq 1$. We also have $a_{bi}\leq a_{bw}$. This gives $$\epsilon^*\leq\frac{\epsilon_{i_0}}{1+a_{bw}+a_{i_0w}+\displaystyle\sum_{\substack{i\in D\\i\neq i_0}}\min\biggl\{1+\epsilon_{i_0},\frac{a_{bw}}{a_{bw}}\biggr\}}=\frac{\epsilon_{i_0}}{1+a_{bw}+a_{i_0w}+n-3}=\frac{\epsilon_{i_0}}{n-2+a_{bw}+a_{i_0w}}.$$
	Now, substituting the value of $\epsilon_{i_0}$ from equation \eqref{4CV1}, we get
	$$\epsilon^*\leq \frac{a_{bw}-a_{bi_0}\times a_{i_0w}}{(a_{bi_0}+2)(n-2+a_{bw}+a_{i_0w})}.$$ Note that $\frac{a_{bw}-a_{bi_0}\times a_{i_0w}}{(a_{bi_0}+2)(n-2+a_{bw}+a_{i_0w})}$ decreases as $a_{bi_0}$ or $a_{i_0w}$ increases. Since $a_{bi_0},a_{i_0w}\geq 1$, we get $$\frac{a_{bw}-a_{bi_0}\times a_{i_0w}}{(a_{bi_0}+2)(n-2+a_{bw}+a_{i_0w})}\leq \frac{a_{bw}-1\times 1}{(1+2)(n-2+a_{bw}+1)}=\frac{a_{bw}-1}{3(n-1+a_{bw})},$$ and thus, $\epsilon^*\leq \frac{a_{bw}-1}{3(n-1+a_{bw})}.$\\\\
	Now, consider the PCS $(A_b',A_w')$. Here, $D_1=\{2\}$, $D_2=\phi$ and $D_3=\{3,4,\ldots,n-1\}$. By equation \eqref{4CV1}, we get $\epsilon_{2}=\frac{a_{bw}-1}{3}$ and $\epsilon_i=0$ for $i=3,4,\ldots,n-1$. So, by equation \eqref{4max}, $\eta=\epsilon_{2}$. Since $2\in D_1$, by statement \textit{1} of Theorem \ref{4exact_obj}, we get $\epsilon^*=\frac{a_{bw}-1}{3(n-1+a_{bw})}$. Hence the proof.
\end{proof}
\hspace{-0.7cm}
Theorem \eqref{4ci} implies that, for the Saaty scale,
\begin{equation}\label{4value_ci}
	\text{CI}_{a_{bw}}(n)=\begin{cases}
		\begin{aligned}
		&\max\biggl\{\frac{a_{bw}-1}{3(2+a_{bw})},\frac{a_{bw}(a_{bw}-1)}{2a_{bw}^2+5a_{bw}+2}\biggr\}
		\end{aligned}\quad \quad\quad \quad \quad \quad\quad \quad\quad  \ \quad\quad\ \text{ if } n=3,\\
		\begin{aligned}
		&\max\biggl\{\frac{a_{bw}-1}{3(n-1+a_{bw})},\frac{a_{bw}(a_{bw}-1)}{2a_{bw}^2+(3n-4)a_{bw}+2}\bigg.,\\
		&\biggl.\quad \ \ \quad\frac{a_{bw}^2-1}{3a_{bw}^2+6a_{bw}+7+(n-4)\min\{a_{bw}^2+a_{bw}+2,1+3a_{bw}\}}\biggr\}
		\end{aligned}\ \ \text{ if } n\geq 4.
	\end{cases}
\end{equation}
Now, $\frac{a_{bw}(a_{bw}-1)}{2a_{bw}^2+(3n-4)a_{bw}+2}=\frac{a_{bw}-1}{2a_{bw}+3n-4+\frac{2}{a_{bw}}}$. Note that for $a_{bw}\geq 2$, $\frac{2}{a_{bw}}\leq 1 < 1+a_{bw}$ holds. This gives ${2a_{bw}+3n-4+\frac{2}{a_{bw}}} < 3n-3+3a_{bw}= 3(n-1+a_{bw})$. So, $\frac{a_{bw}-1}{3(n-1+a_{bw})}<\frac{a_{bw}-1}{2a_{bw}+3n-4+\frac{2}{a_{bw}}}$ for all $n\geq 3$, and thus, we get
\begin{equation}\label{4value_ci_final}
	\text{CI}_{a_{bw}}(n)=\begin{cases}
		\begin{aligned}
			&\frac{a_{bw}(a_{bw}-1)}{2a_{bw}^2+5a_{bw}+2}
		\end{aligned}\quad \quad\quad \quad\quad \quad\quad \ \quad\quad \quad\quad\quad \quad\quad \quad \ \quad\quad \quad\ \quad \ \text{ if } n=3,\\
		\begin{aligned}
			&\max\biggl\{\frac{a_{bw}(a_{bw}-1)}{2a_{bw}^2+(3n-4)a_{bw}+2}\bigg.,\\
			&\biggl.\quad \ \ \quad\frac{a_{bw}^2-1}{3a_{bw}^2+6a_{bw}+7+(n-4)\min\{a_{bw}^2+a_{bw}+2,1+3a_{bw}\}}\biggr\}
		\end{aligned}\ \ \text{ if } n\geq 4.
	\end{cases}
\end{equation}
The values of CI$_{a_{bw}}(n)$ for $a_{bw}=2,3,\ldots,9$ and $n=3,4,\ldots,10$ are provided in Table \ref{4ci_table}.
\begin{table}[H]
	\caption{The values of CI$_{a_{bw}}(n)$}\label{4ci_table}
	\centering		
	\begin{tabular}{@{}cccccccccc@{}}
		\toprule[0.1em]
		$n$&\phantom{}&\multicolumn{8}{c}{$a_{bw}$}\\
		\cmidrule{3-10}
		&&$2$&$3$&$4$&$5$&$6$&$7$&$8$&$9$\\
		\midrule
		$3$&&$0.1$&$0.1714$&$0.2222$&$0.2597$&$0.2885$&$0.3111$&$0.3294$&$0.3445$\\
		$4$&&$0.0968$&$0.1538$&$0.1899$&$0.2174$&$0.2459$&$0.2692$&$0.2887$&$0.3051$\\
		$5$&&$0.0789$&$0.1290$&$0.1630$&$0.1875$&$0.2143$&$0.2373$&$0.2569$&$0.2738$\\
		$6$&&$0.0667$&$0.1111$&$0.1429$&$0.1667$&$0.1899$&$0.2121$&$0.2314$&$0.2483$\\
		$7$&&$0.0577$&$0.0976$&$0.1271$&$0.15$&$0.1705$&$0.1918$&$0.2105$&$0.2271$\\
		$8$&&$0.0508$&$0.0870$&$0.1145$&$0.1364$&$0.1546$&$0.175$&$0.1931$&$0.2093$\\
		$9$&&$0.0455$&$0.0784$&$0.1042$&$0.125$&$0.1423$&$0.1609$&$0.1783$&$0.1941$\\
		$10$&&$0.0411$&$0.0714$&$0.0955$&$0.1154$&$0.1321$&$0.1489$&$0.1657$&$0.1809$\\			
		\bottomrule[0.1em]				
	\end{tabular}
\end{table}
\begin{remark}
	It is important to note that while Theorem \ref{4ci} applies exclusively to the Saaty scale, the same approach can easily be replicated for other scales as well.
\end{remark}
\subsection{Numerical Examples}
In this subsection, we illustrate and validate the proposed approach through numerical examples. For each example, we compute the optimal weights and analyse the sensitivity of the given PCS.\\\\
\textbf{Example 1: }Let $C=\{c_1,c_2,\ldots,c_5\}$ be the set of decision criteria with $c_1$ as the best and $c_5$ as the worst criterion. Let $A_b=(1,2,3,4,7)$ be the best-to-other vector, and $A_w=(7,2,3,2,1)^T$ be the other-to-worst vector.\\\\
\textbf{Calculation of Weights}\\
Here, $D_1=\{2\}$, $D_2=\{3,4\}$ and $D_3=\phi$. So, by equation \eqref{4max}, $\eta=\max\{\epsilon_2,\epsilon_3,\epsilon_4,\epsilon_{2,3},\epsilon_{2,4}\}$. Now, equation \eqref{4CV} gives $\epsilon_2=0.75$, $\epsilon_3=0.4$, $\epsilon_4=0.1667$, $\epsilon_{2,3}=0.7143$ and $\epsilon_{2,4}=0.5$, and thus, $\eta=\epsilon_2=0.75$. Since $2\in D_1$, by statement \textit{1} of Theorem \ref{4exact_obj}, we get
\begin{eqnarray*}
	\sigma&=&1+a_{15}+a_{25}+\min\biggl\{a_{35}+\epsilon_{2},\frac{a_{15}}{a_{13}}\biggr\}+\min\biggl\{a_{45}+\epsilon_{2},\frac{a_{15}}{a_{14}}\biggr\}\\
	&=&1+7+2+\min\biggl\{3+0.75,\frac{7}{3}\biggr\}+\min\biggl\{2+0.75,\frac{7}{4}\biggr\}=14.0833,\\
	w_1&=&\frac{a_{15}-\epsilon_{2}}{\sigma}=\frac{7-0.75}{14.0833}=0.4438,\\
	w_2&=&\frac{a_{25}+\epsilon_{2}}{\sigma}=\frac{2+0.75}{14.0833}=0.1953,\\
	w_3&=&\frac{\min\biggl\{a_{35}+\epsilon_{2},\frac{a_{15}}{a_{13}}\biggr\}}{\sigma}=\frac{\min\biggl\{3+0.75,\frac{7}{3}\biggr\}}{14.0833}=0.1657,\\
	w_4&=&\frac{\min\biggl\{a_{45}+\epsilon_{2},\frac{a_{15}}{a_{14}}\biggr\}}{\sigma}=\frac{\min\biggl\{2+0.75,\frac{7}{4}\biggr\}}{14.0833}=0.1243,\\
	w_5&=&\frac{1}{\sigma}=\frac{1}{14.0833}=0.0710,\\
	\epsilon^*&=&\frac{\epsilon_{2}}{\sigma}=\frac{0.75}{14.0833}=0.0533.
\end{eqnarray*}
Now, from equation \eqref{4cr} and Table \ref{4ci_table}, we get CR $=\frac{0.0533}{0.2373}=0.2246$.\\\\
\textbf{Sensitivity of Data}\\
Statement \textit{1} of Corollary \ref{4sa} establishes that for all four PCSs $A_b=(1,2,3,4,7)$, $A_w=(7,2,a,b,1)^T$, where $a\in\{2,3\}$, $b\in\{1,2\}$, both the optimal weight set and $\epsilon^*$ remain the same.\\\\ 
\textbf{Example 2: }Let $C=\{c_1,c_2,\ldots,c_5\}$ be the set of decision criteria with $c_1$ as the best and $c_5$ as the worst criterion. Let $A_b=(1,4,3,2,9)$ be the best-to-other vector, and $A_w=(9,2,4,7,1)^T$ be the other-to-worst vector.\\\\
\textbf{Calculation of Weights}\\
Here, $D_1=\{2\}$, $D_2=\{3,4\}$ and $D_3=\phi$. So, by equation \eqref{4max}, $\eta=\max\{\epsilon_2,\epsilon_3,\epsilon_4,\epsilon_{2,3},\epsilon_{2,4}\}$. Now, equation \eqref{4CV} gives $\epsilon_2=0.1667$, $\epsilon_3=0.6$, $\epsilon_4=1.25$, $\epsilon_{2,3}=0.4444$ and $\epsilon_{2,4}=0.75$, and thus, $\eta=\epsilon_4=1.25$. Since $4\in D_2$, by statement \textit{2} of Theorem \ref{4exact_obj}, we get
\begin{eqnarray*}
	\sigma&=&1+a_{15}+a_{45}+\min\biggl\{a_{25}+\epsilon_{4},\frac{a_{15}+2\epsilon_{4}}{a_{12}}\biggr\}+\min\biggl\{a_{35}+\epsilon_{4},\frac{a_{15}+2\epsilon_{4}}{a_{13}}\biggr\}\\
	&=&1+9+7+\min\biggl\{2+1.25,\frac{9+2\times 1.25}{4}\biggr\}+\min\biggl\{4+1.25,\frac{9+2\times 1.25}{3}\biggr\}\\
	&=&23.7083,\\
	w_1&=&\frac{a_{15}+\epsilon_{4}}{\sigma}=\frac{9+1.25}{23.7083}=0.4323,\\
	w_2&=&\frac{\min\biggl\{a_{25}+\epsilon_{4},\frac{a_{15}+2\epsilon_{4}}{a_{12}}\biggr\}}{\sigma}=\frac{\min\biggl\{2+1.25,\frac{9+2\times 1.25}{4}\biggr\}}{23.7083}=0.1213,\\
	w_3&=&\frac{\min\biggl\{a_{35}+\epsilon_{4},\frac{a_{15}+2\epsilon_{4}}{a_{13}}\biggr\}}{\sigma}=\frac{\min\biggl\{4+1.25,\frac{9+2\times 1.25}{3}\biggr\}}{23.7083}=0.1617,\\
	w_4&=&\frac{a_{45}-\epsilon_{4}}{\sigma}=\frac{7-1.25}{23.7083}=0.2425,\\
	w_5&=&\frac{1}{\sigma}=\frac{1}{23.7083}=0.0422,\\
	\epsilon^*&=&\frac{\epsilon_{4}}{\sigma}=\frac{1.25}{23.7083}=0.0527.
\end{eqnarray*}
Now, from equation \eqref{4cr} and Table \ref{4ci_table}, we get CR $=\frac{0.0527}{0.2738}=0.1925$.\\\\
\textbf{Sensitivity of Data}\\ 
Statement \textit{2} of Corollary \ref{4sa} establishes that for all nine PCSs $A_b=(1,4,3,2,9)$, $A_w=(9,a,b,7,1)^T$, where $a\in\{2,3,4\}$, $b\in\{3,4,5\}$, both the optimal weight set and $\epsilon^*$ remain the same.\\\\
\textbf{Example 3: }Let $C=\{c_1,c_2,\ldots,c_7\}$ be the set of decision criteria with $c_1$ as the best and $c_7$ as the worst criterion. Let $A_b=(1,1,4,3,2,4,5)$ be the best-to-other vector, and $A_w=(5,2,5,2,3,2,1)^T$ be the other-to-worst vector.\\\\
\textbf{Calculation of Weights}\\
Here, $D_1=\{2\}$, $D_2=\{3,4,5,6\}$ and $D_3=\phi$. So, by equation \eqref{4max}, $\eta=\max \left\{\epsilon_2,\epsilon_3,\epsilon_4,\epsilon_5,\epsilon_6,\epsilon_{2,3},\right.$ $\left.\epsilon_{2,4},\epsilon_{2,5},\epsilon_{2,6}\right\}$. Now, equation \eqref{4CV} gives $\epsilon_2=1$, $\epsilon_3=2.5$, $\epsilon_4=0.2$, $\epsilon_5=0.25$, $\epsilon_6=0.8333$, $\epsilon_{2,3}=2.5714$, $\epsilon_{2,4}=0.6667$, $\epsilon_{2,5}=0.8$ and $\epsilon_{2,6}=0.8571$, and thus, $\eta=\epsilon_{2,3}=2.5714$. By statement \textit{3} of Theorem \ref{4exact_obj}, we get
\begin{eqnarray*}
	\sigma&=&1+a_{27}+a_{37}+a_{12}\times a_{27}+(a_{12}+1)\times \epsilon_{2,3}\\
	&&\quad\quad\quad\quad+\min\biggl\{a_{47}+\epsilon_{2,3},\frac{a_{12}\times a_{27}+(a_{12}+2)\times \epsilon_{2,3}}{a_{14}}\biggr\}\\
	&&\quad\quad\quad\quad\quad\quad\quad\quad+\min\biggl\{a_{57}+\epsilon_{2,3},\frac{a_{12}\times a_{27}+(a_{12}+2)\times \epsilon_{2,3}}{a_{15}}\biggr\}\\
	&&\quad\quad\quad\quad\quad\quad\quad\quad\quad\quad\quad\quad+\min\biggl\{a_{67}+\epsilon_{2,3},\frac{a_{12}\times a_{27}+(a_{12}+2)\times \epsilon_{2,3}}{a_{16}}\biggr\}\\
	&=&1+2+5+1\times 2+(1+1)\times 2.5714\\
	&&\quad\quad\quad\quad+\min\biggl\{2+2.5714,\frac{1\times 2+(1+2)\times 2.5714}{3}\biggr\}\\
	&&\quad\quad\quad\quad\quad\quad\quad\quad+\min\biggl\{3+2.5714,\frac{1\times 2+(1+2)\times 2.5714}{2}\biggr\}\\
	&&\quad\quad\quad\quad\quad\quad\quad\quad\quad\quad\quad\quad+\min\biggl\{2+2.5714,\frac{1\times 2+(1+2)\times 2.5714}{4}\biggr\}\\
	&=&25.6667,\\
	w_1&=&\frac{a_{12}\times a_{27}+(a_{12}+1)\times \epsilon_{2,3}}{\sigma}=\frac{1\times 2+(1+1)\times 2.5714}{25.6667}=0.2783,\\
	w_2&=&\frac{a_{27}+\epsilon_{2,3}}{\sigma}=\frac{2+2.5714}{25.6667}=0.1781,\\
	w_3&=&\frac{a_{37}-\epsilon_{2,3}}{\sigma}=\frac{5-2.5714}{25.6667}=0.0946,\\
	w_4&=&\frac{\min\biggl\{a_{47}+\epsilon_{2,3},\frac{a_{12}\times a_{27}+(a_{12}+2)\times \epsilon_{2,3}}{a_{14}}\biggr\}}{\sigma}=\frac{\min\biggl\{2+2.5714,\frac{1\times 2+(1+2)\times 2.5714}{3}\biggr\}}{25.6667}=0.1262,\\
	w_5&=&\frac{\min\biggl\{a_{57}+\epsilon_{2,3},\frac{a_{12}\times a_{27}+(a_{12}+2)\times \epsilon_{2,3}}{a_{15}}\biggr\}}{\sigma}=\frac{\min\biggl\{3+2.5714,\frac{1\times 2+(1+2)\times 2.5714}{2}\biggr\}}{25.6667}=0.1892,\\
	w_6&=&\frac{\min\biggl\{a_{67}+\epsilon_{2,3},\frac{a_{12}\times a_{27}+(a_{12}+2)\times \epsilon_{2,3}}{a_{16}}\biggr\}}{\sigma}=\frac{\min\biggl\{2+2.5714,\frac{1\times 2+(1+2)\times 2.5714}{4}\biggr\}}{25.6667}=0.0946,\\
	w_7&=&\frac{1}{\sigma}=\frac{1}{25.6667}=0.0390,\\
	\epsilon^*&=&\frac{\epsilon_{2,3}}{\sigma}=\frac{2.5714}{25.6667}=0.1002.
\end{eqnarray*}
Now, from equation \eqref{4cr} and Table \ref{4ci_table}, we get CR $=\frac{0.1002}{0.15}=0.6680$.\\\\
\textbf{Sensitivity of Data}\\
Statement \textit{3} of Corollary \ref{4sa} establishes that for all seventy-five PCSs $A_b=(1,1,4,3,2,$ $4,5)$, $A_w=(5,2,5,a,b,d,1)^T$, where $a\in\{1,2,\ldots,5\}$, $b\in\{3,4,5\}$, $d=\{1,2,\ldots,5\}$, both the optimal weight set and $\epsilon^*$ remain the same.\\\\
\textbf{Example 4: }Let $C=\{c_1,c_2,c_3,c_4\}$ be the set of decision criteria with $c_1$ as the best and $c_4$ as the worst criterion. Let $A_b=(1,5,4,8)$ be the best-to-other vector, and $A_w=(8,4,1,1)^T$ be the other-to-worst vector.\\\\
\textbf{Calculation of Weights}\\
Here, $D_1=\{3\}$, $D_2=\{2\}$ and $D_3=\phi$. So, by equation \eqref{4max}, $\eta=\max\{\epsilon_2,\epsilon_3,\epsilon_{3,2}\}$. Now, equation \eqref{4CV} gives $\epsilon_2=1.7143$, $\epsilon_3=0.6667$ and $\epsilon_{3,2}=1.4545$, and thus, $\eta=\epsilon_2=1.7143$. Since $2 \in D_2$, by statement \textit{2} of Theorem \ref{4exact_obj}, we get
\begin{eqnarray*}
	\sigma&=&1+a_{14}+a_{24}+\min\biggl\{a_{34}+\epsilon_{2},\frac{a_{14}+2\epsilon_{2}}{a_{13}}\biggr\}\\
	&=&1+8+4+\min\biggl\{1+1.7143,\frac{8+2\times 1.7143}{4}\biggr\}\\
	&=&15.7143,\\
	w_1&=&\frac{a_{14}+\epsilon_{2}}{\sigma}=\frac{8+1.7143}{15.7143}=0.6182,\\
	w_2&=&\frac{a_{24}-\epsilon_{2}}{\sigma}=\frac{4-1.7143}{15.7143}=0.1455,\\
	w_3&=&\frac{\min\biggl\{a_{34}+\epsilon_{2},\frac{a_{14}+2\epsilon_{2}}{a_{13}}\biggr\}}{\sigma}=\frac{\min\biggl\{1+1.7143,\frac{8+2\times 1.7143}{4}\biggr\}}{15.7143}=0.1727,\\
	w_4&=&\frac{1}{\sigma}=\frac{1}{15.7143}=0.0636,\\
	\epsilon^*&=&\frac{\epsilon_{2}}{\sigma}=\frac{1.7143}{15.7143}=0.1091.
\end{eqnarray*}
Now, from equation \eqref{4cr} and Table \ref{4ci_table}, we get CR $=\frac{0.1091}{0.2887}=0.3779$.\\\\
\textbf{Sensitivity of Data}\\
Statement \textit{2} of Corollary \ref{4sa} establishes that for both PCSs $A_b=(1,5,a,8)$, $A_w=(8,4,1,1)^T$, where $a\in\{3,4\}$, both the optimal weight set and $\epsilon^*$ remain the same.\\\\
\textbf{Example 5: }Let $C=\{c_1,c_2,\ldots,c_5\}$ be the set of decision criteria with $c_1$ as the best and $c_5$ as the worst criterion. Let $A_b=(1,6,3,4,6)$ be the best-to-other vector, and $A_w=(6,6,2,1,1)^T$ be the other-to-worst vector.\\\\
\textbf{Calculation of Weights}\\
Here, $D_1=\{4\}$, $D_2=\{2\}$ and $D_3=\{3\}$. So, by equation \eqref{4max}, $\eta=\max\{\epsilon_2,\epsilon_4,\epsilon_{4,2}\}$. Now, equation \eqref{4CV} gives $\epsilon_2=3.75$, $\epsilon_4=0.3333$ and $\epsilon_{4,2}=2.6667$, and thus, $\eta=\epsilon_2=3.75$. Since $2\in D_2$, by statement \textit{2} of Theorem \ref{4exact_obj}, we get
\begin{eqnarray*}
	\sigma&=&1+a_{15}+a_{25}+\min\biggl\{a_{35}+\epsilon_{2},\frac{a_{15}+2\epsilon_{2}}{a_{13}}\biggr\}+\min\biggl\{a_{45}+\epsilon_{2},\frac{a_{15}+2\epsilon_{2}}{a_{14}}\biggr\}\\
	&=&1+6+6+\min\biggl\{2+3.75,\frac{6+2\times 3.75}{3}\biggr\}+\min\biggl\{1+3.75,\frac{6+2\times 3.75}{4}\biggr\}\\
	&=&20.875,\\
	w_1&=&\frac{a_{15}+\epsilon_{2}}{\sigma}=\frac{6+3.75}{20.875}=0.4671,\\
	w_2&=&\frac{a_{25}-\epsilon_{2}}{\sigma}=\frac{6-3.75}{20.875}=0.1078,\\
	w_3&=&\frac{\min\biggl\{a_{35}+\epsilon_{2},\frac{a_{15}+2\epsilon_{2}}{a_{13}}\biggr\}}{\sigma}=\frac{\min\biggl\{2+3.75,\frac{6+2\times 3.75}{3}\biggr\}}{20.875}=0.2156,\\
	w_4&=&\frac{\min\biggl\{a_{45}+\epsilon_{2},\frac{a_{15}+2\epsilon_{2}}{a_{14}}\biggr\}}{\sigma}=\frac{\min\biggl\{1+3.75,\frac{6+2\times 3.75}{4}\biggr\}}{20.875}=0.1617,\\
	w_5&=&\frac{1}{\sigma}=\frac{1}{20.875}=0.0479,\\
	\epsilon^*&=&\frac{\epsilon_{2}}{\sigma}=\frac{3.75}{20.875}=0.1796.
\end{eqnarray*}
Now, from equation \eqref{4cr} and Table \ref{4ci_table}, we get CR $=\frac{0.1796}{0.2143}=0.8381$.\\\\
\textbf{Sensitivity of Data}\\
Statement \textit{2} of Corollary \ref{4sa} establishes that for all thirty-six PCSs $A_b=(1,6,3,4,6)$, $A_w=(6,6,a,b,1)^T$, where $a\in\{1,2,\ldots,6\}$, $b\in\{1,2,\ldots,6\}$, both the optimal weight set and $\epsilon^*$ remain the same.\\
\begin{remark}
	The most commonly used optimization tool for the linear BWM is the Excel solver provided on \url{https://bestworstmethod.com/software/}. However, the optimal weights obtained using this solver ($w_1^*=0.4706$, $w_2^*=0.1176$, $w_3^*=0.2353$, $w_4^*=0.0588$ and $w_5^*=0.1176$, with $\epsilon^*=0.5882$) differ from the actual optimal weights. This discrepancy raises concerns about the reliability of the solver.
\end{remark}
\section{A Real-World Application}
In this section, we demonstrate the applicability of the proposed approach by ranking the drivers of Industry 4.0, sustainability and Circular Economy (CE) in relation to Sustainable Development Goals (SDGs)-driven Agri-Food Supply Chains (AFSCs).\\\\
The AFSC represents the complete journey of agricultural products from farm to fork, encompassing production, processing, packaging, distribution, retail and waste management. This complex system plays a crucial role in global food security, economic development and environmental sustainability \cite{fao2021}.\\\\
The SDGs are 17 global goals set by the United Nations (2015–2030) to address poverty, inequality, climate change and sustainable development \cite{un2015}. Several SDGs are deeply interconnected with AFSC operations, with multiple goals directly addressing critical aspects of sustainable food systems. SDG 1 (no poverty) strengthens rural livelihoods by supporting smallholder farmers and fair agricultural value chains. SDG 2 (zero hunger) serves as a foundational goal, driving efforts to ensure food security through sustainable agricultural practices. SDG 3 (good health and well-being) safeguards nutrition security and food safety throughout the supply chain. Environmental sustainability is addressed through SDG 6 (clean water and sanitation), which promotes efficient water management in agriculture, and SDG 7 (affordable and clean energy), which encourages renewable energy adoption in food processing and transportation. SDG 12 (responsible consumption and production) is particularly crucial as it fosters CE principles to minimize food waste and optimize resource use. Climate considerations are embedded in SDG 13 (climate action), which targets emission reductions across agricultural activities, and SDG 15 (life on land), which protects vital ecosystems from unsustainable farming practices. This intrinsic connection between AFSCs and SDGs highlights how optimizing food supply systems can simultaneously advance multiple global development objectives while creating more resilient and equitable food networks.\\\\
Industry 4.0 technologies are revolutionizing AFSCs by making them smarter and more efficient. IoT sensors enable real-time crop monitoring through continuous data collection on soil conditions and plant health, while AI analytics detect early signs of disease by recognizing patterns in sensor and image data \cite{wolfert2017big,kayikci2022food}. Blockchain technology ensures end-to-end traceability by maintaining secure, immutable records of production and distribution processes \cite{sharma2024developing}. These innovations directly contribute to SDG 2 through yield optimization and SDG 12 by minimizing food losses. Furthermore, automated irrigation systems and precision agriculture tools promote SDG 6 through reduced water usage, while smart energy management in processing and transportation supports SDG 7 \cite{alam2023analysis}. These technological innovations are not merely improvements but vital necessities for developing AFSCs capable of meeting the dual challenges of global food security and environmental sustainability, making implementation of Industry 4.0 crucial for realizing the 2030 SDG agenda in the agriculture sector.\\\\
Effective waste management is crucial in AFSCs to mitigate environmental degradation and economic inefficiencies. The CE serves as a transformative approach to minimize waste and maximize resource efficiency throughout food value chains. As a critical element of sustainable AFSCs, CE enhances resource utilization to strengthen food security, directly supporting SDG 2 \cite{perccin2025evaluating}. Through innovative waste valorization and reusable packaging solutions, CE contributes significantly to SDG 12 \cite{barros2020mapping}. By converting agricultural residues into bioenergy, CE reduces greenhouse gas emissions, tangibly furthering SDG 13 \cite{zhang2022circular}. By converting food byproducts into valuable materials like biofuels and organic fertilizers, and implementing sustainable packaging systems, CE addresses critical challenges in AFSCs - reducing post-harvest losses establishing closed-loop systems that support SDG 15, and creating new economic opportunities for smallholder farmers, thereby advancing SDG 1 \cite{barros2020mapping}. Consequently, CE principles are no longer optional but imperative for AFSCs to achieve 2030 SDG targets while fostering resilient, future-proof food systems. \\\\
In this study, we rank eighteen drivers from three categories - Industry 4.0, sustainability, and CE (listed in Table \ref{4drivers}) - based on their criticality for achieving SDGs in AFSCs. The pairwise comparisons between categories and among drivers within each category, as provided by five homogeneous experts $E_1-E_5$, are adopted from the work of Per{\c{c}}in et al. \cite{perccin2025evaluating}. Based on these comparisons, the category weights and the local weights for each driver are calculated for each expert. Using these weights, the global weights for each driver are calculated for each expert. The final weights of drivers are determined by aggregating the global weights across all experts using the arithmetic mean method, enabling the ranking of all drivers. The category weights are presented in Table \ref{4weights_driver}, while the local weights and the global weights for each expert, the final weights, and the ranking are given in Table \ref{4weights_subdriver}.\\\\
Our analysis reveals that ``resource efficiency" is the most critical and ``artificial intelligence" is the least critical driver among the evaluated drivers. This ranking provides valuable insights for prioritizing implementation efforts in AFSCs to effectively advance SDG achievement.
\begin{table}[H]
	\caption{List of drivers assessed for SDG impact in AFSCs \cite{perccin2025evaluating}}\label{4drivers}
	\centering		
	\begin{tabular}{@{}ll@{}}
		\toprule[0.1em]
		Category& Driver\\
		\midrule
		$c_1$: Industry 4.0&$c_{11}$: Cloud computing\\
		&$c_{12}$: Sensors and robotics\\
		&$c_{13}$: Internet of Things (IoT) \\
		&$c_{14}$: Blockchain\\
		&$c_{15}$: Big Data Analytics (BDA)\\
		&$c_{16}$: Artificial Intelligence (AI)\\
		$c_2$: Sustainability&$c_{21}$: Economic sustainability\\
		&$c_{22}$: Environmental sustainability\\
		&$c_{23}$: Social sustainability\\
		&$c_{24}$: Innovative business models\\
		&$c_{25}$: Competitiveness\\
		&$c_{26}$: Achievement of standards and SDGs\\
		$c_3$: CE&$c_{31}$: Resource efficiency\\
		&$c_{32}$: Waste and emissions reduction\\
		&$c_{33}$: Supply chain connectivity\\
		&$c_{34}$: Traceability and transparency\\
		&$c_{35}$: Legal compliance\\
		&$c_{36}$: Stakeholders' rights\\
		\bottomrule[0.1em]				
	\end{tabular}
\end{table}
\begin{table}[H]
	\caption{Expert-wise category weights}\label{4weights_driver}
	\centering		
	\begin{tabular}{@{}ccccccc@{}}
		\toprule[0.1em]
		Category&\phantom{}&\multicolumn{5}{c}{Weight}\\
		\cline{3-7}
		&&$E_1$&$E_2$&$E_3$&$E_4$&$E_5$\\
		\midrule
		$c_1$&&$0.1111$&$0.1868$&$0.0667$&$0.0769$&$0.0769$\\
		$c_2$&&$0.2444$&$0.0769$&$0.2533$&$0.7582$&$0.7949$\\
		$c_3$&&$0.6444$&$0.7363$&$0.6800$&$0.1648$&$0.1282$\\
		$\epsilon^*$&&$0.0889$&$0.1978$&$0.0800$&$0.0659$&$0.0949$\\
		CR&&$0.3423$&$0.6358$&$0.2322$&$0.1914$&$0.3823$\\
		\bottomrule[0.1em]				
	\end{tabular}
\end{table}
\begin{sidewaystable}[h]
\caption{Final aggregated weights and ranking of drivers}\label{4weights_subdriver}
\centering		
\begin{tabular}{@{}cccccccccccccccc@{}}
	\toprule[0.1em]
	Driver&&\multicolumn{5}{c}{Local weight}&&\multicolumn{5}{c}{Global weight}&&Final weight&Rank\\
	\cline{3-7}\cline{9-13}
	&&$E_1$&$E_2$&$E_3$&$E_4$&$E_5$&&$E_1$&$E_2$&$E_3$&$E_4$&$E_5$&&&\\
	\midrule
	$c_{11}$&&$0.4900$&$0.0749$&$0.1032$&$0.5164$&$0.1092$&&$0.0544$&$0.0140$&$0.0069$&$0.0397$&$0.0084$&&$0.0247$&$13$\\
	$c_{12}$&&$0.0847$&$0.4329$&$0.1032$&$0.1284$&$0.1092$&&$0.0094$&$0.0809$&$0.0069$&$0.0099$&$0.0084$&&$0.0231$&$14$\\
	$c_{13}$&&$0.1186$&$0.1747$&$0.1721$&$0.0917$&$0.0396$&&$0.0132$&$0.0326$&$0.0115$&$0.0071$&$0.0030$&&$0.0135$&$16$\\
	$c_{14}$&&$0.1977$&$0.0380$&$0.0425$&$0.0917$&$0.1820$&&$0.0220$&$0.0071$&$0.0028$&$0.0071$&$0.0140$&&$0.0106$&$17$\\
	$c_{15}$&&$0.0659$&$0.1747$&$0.4069$&$0.1284$&$0.4509$&&$0.0073$&$0.0326$&$0.0271$&$0.0099$&$0.0347$&&$0.0223$&$15$\\
	$c_{16}$&&$0.0430$&$0.1048$&$0.1721$&$0.0434$&$0.1092$&&$0.0048$&$0.0196$&$0.0115$&$0.0033$&$0.0084$&&$0.0095$&$18$\\
	$\epsilon^*$&&$0.1032$&$0.0911$&$0.1093$&$0.1255$&$0.0949$&&-&-&-&-&-&&-&-\\
	CR&&$0.4155$&$0.3671$&$0.5154$&$0.5054$&$0.3823$&&-&-&-&-&-&&-&-\\
	\midrule
	$c_{21}$&&$0.4280$&$0.1309$&$0.0824$&$0.4209$&$0.1654$&&$0.1046$&$0.0101$&$0.0209$&$0.3191$&$0.1315$&&$0.1172$&$2$\\
	$c_{22}$&&$0.0611$&$0.4800$&$0.1154$&$0.0540$&$0.1654$&&$0.0149$&$0.0369$&$0.0292$&$0.0409$&$0.1315$&&$0.0507$&$9$\\
	$c_{23}$&&$0.1101$&$0.0935$&$0.4487$&$0.1619$&$0.1654$&&$0.0269$&$0.0072$&$0.1137$&$0.1227$&$0.1315$&&$0.0804$&$5$\\
	$c_{24}$&&$0.0340$&$0.1309$&$0.1154$&$0.1619$&$0.3993$&&$0.0083$&$0.0101$&$0.0292$&$0.1227$&$0.3174$&&$0.0975$&$4$\\
	$c_{25}$&&$0.1834$&$0.0339$&$0.1923$&$0.0396$&$0.0336$&&$0.0448$&$0.0026$&$0.0487$&$0.0300$&$0.0267$&&$0.0306$&$12$\\
	$c_{26}$&&$0.1834$&$0.1309$&$0.0458$&$0.1619$&$0.0709$&&$0.0448$&$0.0101$&$0.0116$&$0.1227$&$0.0564$&&$0.0491$&$10$\\
	$\epsilon^*$&&$0.1223$&$0.1745$&$0.1282$&$0.0647$&$0.0970$&&-&-&-&-&-&&-&-\\
	CR&&$0.4925$&$0.7028$&$0.6045$&$0.2608$&$0.3907$&&-&-&-&-&-&&-&-\\
	\midrule
	$c_{31}$&&$0.5442$&$0.1705$&$0.1166$&$0.4680$&$0.5610$&&$0.3507$&$0.1255$&$0.0793$&$0.0771$&$0.0719$&&$0.1409$&$1$\\
	$c_{32}$&&$0.0443$&$0.1705$&$0.1166$&$0.0368$&$0.1295$&&$0.0286$&$0.1255$&$0.0793$&$0.0061$&$0.0166$&&$0.0512$&$8$\\
	$c_{33}$&&$0.1379$&$0.0568$&$0.1943$&$0.1209$&$0.0527$&&$0.0889$&$0.0418$&$0.1321$&$0.0199$&$0.0068$&&$0.0579$&$7$\\
	$c_{34}$&&$0.0766$&$0.3977$&$0.0360$&$0.0864$&$0.0719$&&$0.0494$&$0.2928$&$0.0245$&$0.0142$&$0.0092$&&$0.0780$&$6$\\
	$c_{35}$&&$0.0985$&$0.1023$&$0.4533$&$0.2016$&$0.0925$&&$0.0635$&$0.0753$&$0.3083$&$0.0332$&$0.0119$&&$0.0984$&$3$\\
	$c_{36}$&&$0.0985$&$0.1023$&$0.0833$&$0.0864$&$0.0925$&&$0.0635$&$0.0753$&$0.0566$&$0.0142$&$0.0119$&&$0.0443$&$11$\\
	$\epsilon^*$&&$0.1451$&$0.1136$&$0.1295$&$0.1367$&$0.0863$&&-&-&-&-&-&&-&-\\
	CR&&$0.5845$&$0.6817$&$0.5216$&$0.5506$&$0.3476$&&-&-&-&-&-&&-&-\\
	\bottomrule[0.1em]				
\end{tabular}
\end{sidewaystable}
\section{Conclusions and Future Directions}
The BWM is an emerging MCDM method that researchers have widely applied in various real-world scenarios. In this work, we focus on the linear BWM, one of the most commonly used models of BWM. We propose a framework to derive an analytical expression for the optimal weights, eliminating dependency on optimization software and reducing computational complexity. Beyond these computational advantages, the proposed analytical approach also reveals a critical limitation of the model—its low sensitivity to data variations. Using this framework, we also calculate the consistency index, which facilitate the determination of the consistency ratio—a key indicator of inconsistency in pairwise comparisons. Furthermore, we rank eighteen drivers across three categories - Industry 4.0, sustainability and circular economy - according to their criticality for sustainable development goals-driven agri-food supply chains.\\\\
Within the broader scope of this study, several problems remain unaddressed. For some models of BWM, such as the Euclidean BWM \cite{kocak2018euclidean} and the $\alpha$-FBWM \cite{ratandhara2024alpha}, an analytical expression for optimal weights has yet to be derived. Another key future direction involves determining the threshold values for the consistency ratio in the linear BWM to assess the admissibility of input data, analogous to the work of Liang et al. \cite{liang2020consistency} in case of the nonlinear BWM.
\section*{Acknowledgements}
The first author gratefully acknowledges the Council of Scientific \& Industrial Research (CSIR), India for financial support to carry out the research work.
\section*{Declaration of Conflict of Interest}
The authors declare that they have no known conflict of financial interests or personal relationships that could have appeared to influence the work reported in this paper.
\bibliographystyle{plain}

\begin{thebibliography}{10}
	
	\bibitem{ahmadi2017assessing}
	Ahmadi, H.~B., Kusi-Sarpong, S., and Rezaei, J. (2017).
	\newblock Assessing the social sustainability of supply chains using best worst
	method.
	\newblock {\em Resources, Conservation and Recycling}, 126:99--106. https://doi.org/10.1016/j.resconrec.2017.07.020
	
	\bibitem{alam2023analysis}
	Alam, M. F.~B., Tushar, S.~R., Zaman, S.~M., Gonzalez, E. D.~S., Bari, A.~M.,
	and Karmaker, C.~L. (2023).
	\newblock Analysis of the drivers of agriculture 4.0 implementation in the
	emerging economies: Implications towards sustainability and food security.
	\newblock {\em Green Technologies and Sustainability}, 1(2):100021. https://doi.org/10.1016/j.grets.2023.100021
	
	\bibitem{ali2019hesitant}
	Ali, A., and Rashid, T. (2019).
	\newblock Hesitant fuzzy best-worst multi-criteria decision-making method and
	its applications.
	\newblock {\em International Journal of Intelligent Systems}, 34(8):1953--1967. https://doi.org/10.1002/int.22131
	
	\bibitem{amiri2020goal}
	Amiri, M., and Emamat, M. S. M.~M. (2020).
	\newblock A goal programming model for bwm.
	\newblock {\em Informatica}, 31(1):21--34. https://doi.org/10.15388/20-INFOR389
	
	\bibitem{barros2020mapping}
	Barros, M.~V., Salvador, R., De~Francisco, A.~C., and Piekarski, C.~M. (2020).
	\newblock Mapping of research lines on circular economy practices in
	agriculture: From waste to energy.
	\newblock {\em Renewable and Sustainable Energy Reviews}, 131:109958. https://doi.org/10.1016/j.rser.2020.109958
	
	\bibitem{brunelli2019multiplicative}
	Brunelli, M., and Rezaei, J. (2019).
	\newblock A multiplicative best--worst method for multi-criteria decision
	making.
	\newblock {\em Operations Research Letters}, 47(1):12--15. https://doi.org/10.1016/j.orl.2018.11.008
	
	\bibitem{corrente2024better}
	Corrente, S., Greco, S., and Rezaei, J. (2024).
	\newblock Better decisions with less cognitive load: The parsimonious bwm.
	\newblock {\em Omega}, 126:103075. https://doi.org/10.1016/j.omega.2024.103075
	
	
	\bibitem{dawood2023novel}
	Dawood, K.~A., Zaidan, A., Sharif, K.~Y., Ghani, A.~A., Zulzalil, H., and
	Zaidan, B. (2023).
	\newblock Novel multi-perspective usability evaluation framework for selection
	of open source software based on bwm and group vikor techniques.
	\newblock {\em International Journal of Information Technology \& Decision
		Making}, 22(01):187--277. https://doi.org/10.1142/S0219622021500139
	
	\bibitem{fao2021}
	Food and Agriculture Organization (FAO) (2021).
	\newblock Sustainable food systems: concept and framework. https://openknowledge.fao.org/server/api/core/bitstreams/b620989c-407b-4caf-a152-f790f55fec71/content
	
	\bibitem{guo2017fuzzy}
	Guo, S., and Zhao, H. (2017).
	\newblock Fuzzy best-worst multi-criteria decision-making method and its
	applications.
	\newblock {\em Knowledge-Based Systems}, 121:23--31. https://doi.org/10.1016/j.knosys.2017.01.010
	
	\bibitem{gupta2017developing}
	Gupta, P., Anand, S., and Gupta, H. (2017).
	\newblock Developing a roadmap to overcome barriers to energy efficiency in
	buildings using best worst method.
	\newblock {\em Sustainable Cities and Society}, 31:244--259. https://doi.org/10.1016/j.scs.2017.02.005
	
	\bibitem{kayikci2022food}
	Kayikci, Y., Subramanian, N., Dora, M., and Bhatia, M.~S. (2022).
	\newblock Food supply chain in the era of industry 4.0: Blockchain technology
	implementation opportunities and impediments from the perspective of people,
	process, performance, and technology.
	\newblock {\em Production planning \& control}, 33(2-3):301--321. https://doi.org/10.1080/09537287.2020.1810757
	
	\bibitem{kheybari2020sustainable}
	Kheybari, S., Davoodi~Monfared, M., Farazmand, H., and Rezaei, J. (2020).
	\newblock Sustainable location selection of data centers: developing a
	multi-criteria set-covering decision-making methodology.
	\newblock {\em International Journal of Information Technology \& Decision
		Making}, 19(03):741--773. https://doi.org/10.1142/S0219622020500157
	
	\bibitem{kocak2018euclidean}
	Kocak, H., Caglar, A., and Oztas, G.~Z. (2018).
	\newblock Euclidean best--worst method and its application.
	\newblock {\em International Journal of Information Technology \& Decision
		Making}, 17(05):1587--1605. https://doi.org/10.1142/S0219622018500323
	
	\bibitem{lei2022preference}
	Lei, Q., Wu, G., and Wu, Z. (2022).
	\newblock Preference rationality analysis for the best--worst method and its
	application to quality assessment.
	\newblock {\em Computers \& Industrial Engineering}, 174:108758. https://doi.org/10.1016/j.cie.2022.108758
	
	\bibitem{liang2020consistency}
	Liang, F., Brunelli, M., and Rezaei, J. (2020).
	\newblock Consistency issues in the best worst method: Measurements and
	thresholds.
	\newblock {\em Omega}, 96:102175. https://doi.org/10.1016/j.omega.2019.102175
	
	\bibitem{liang2022best}
	Liang, F., Brunelli, M., and Rezaei, J. (2022).
	\newblock Best-worst tradeoff method.
	\newblock {\em Information Sciences}, 610:957--976. https://doi.org/10.1016/j.ins.2022.07.097
	
	\bibitem{mohammadi2020bayesian}
	Mohammadi, M., and Rezaei, J. (2020).
	\newblock Bayesian best-worst method: A probabilistic group decision making
	model.
	\newblock {\em Omega}, 96:102075. https://doi.org/10.1016/j.omega.2019.06.001
	
	\bibitem{mou2016intuitionistic}
	Mou, Q., Xu, Z., and Liao, H. (2016).
	\newblock An intuitionistic fuzzy multiplicative best-worst method for
	multi-criteria group decision making.
	\newblock {\em Information Sciences}, 374:224--239. https://doi.org/10.1016/j.ins.2016.08.074
	
	\bibitem{perccin2025evaluating}
	Per{\c{c}}in, S., Bayraktar, Y., and Kumar, V. (2025).
	\newblock Evaluating the adoption of industry 4.0, sustainability and circular
	economy drivers to achieve sustainable development goals--oriented agri-food
	supply chains.
	\newblock {\em Business Strategy and the Environment}, 34(4):4871--4895. https://doi.org/10.1002/bse.4221
	
	\bibitem{ratandhara2024alpha}
	Ratandhara, H.~M., and Kumar, M. (2024).
	\newblock An $\alpha$-cut intervals based fuzzy best--worst method for
	multi-criteria decision-making.
	\newblock {\em Applied Soft Computing}, 159:111625. https://doi.org/10.1016/j.asoc.2024.111625
	
	\bibitem{ratandhara2024analytical}
	Ratandhara, H.~M., and Kumar, M. (2024b).
	\newblock An analytical framework for the multiplicative best-worst method.
	\newblock {\em Journal of Multi-Criteria Decision Analysis}, 31(5-6):e1840. https://doi.org/10.1002/mcda.1840
	
	\bibitem{rezaei2015best}
	Rezaei, J. (2015).
	\newblock Best-worst multi-criteria decision-making method.
	\newblock {\em Omega}, 53:49--57. https://doi.org/10.1016/j.omega.2014.11.009
	
	\bibitem{rezaei2016best}
	Rezaei, J. (2016).
	\newblock Best-worst multi-criteria decision-making method: Some properties and
	a linear model.
	\newblock {\em Omega}, 64:126--130. https://doi.org/10.1016/j.omega.2015.12.001
	
	\bibitem{rezaei2020concentration}
	Rezaei, J. (2020).
	\newblock A concentration ratio for nonlinear best worst method.
	\newblock {\em International Journal of Information Technology \& Decision
		Making}, 19(03):891--907. https://doi.org/10.1142/S0219622020500170
	
	\bibitem{rezaei2018quality}
	Rezaei, J., Kothadiya, O., Tavasszy, L., and Kroesen, M. (2018).
	\newblock Quality assessment of airline baggage handling systems using servqual
	and bwm.
	\newblock {\em Tourism Management}, 66:85--93. https://doi.org/10.1016/j.tourman.2017.11.009
	
	\bibitem{saaty1994make}
	Saaty, T.~L. (1994).
	\newblock How to make a decision: the analytic hierarchy process.
	\newblock {\em Interfaces}, 24(6):19--43. https://doi.org/10.1287/inte.24.6.19
	
	\bibitem{saaty2004decision}
	Saaty, T.~L. (2004).
	\newblock Decision making—the analytic hierarchy and network processes
	(ahp/anp).
	\newblock {\em Journal of systems science and systems engineering}, 13:1--35. https://doi.org/10.1007/s11518-006-0151-5
	
	\bibitem{safarzadeh2018group}
	Safarzadeh, S., Khansefid, S., and Rasti-Barzoki, M. (2018).
	\newblock A group multi-criteria decision-making based on best-worst method.
	\newblock {\em Computers \& Industrial Engineering}, 126:111--121. https://doi.org/10.1016/j.cie.2018.09.011
		
	\bibitem{sharma2024developing}
	Sharma, A., Bhatia, T., Singh, R.~K., and Sharma, A. (2024).
	\newblock Developing the framework of blockchain-enabled agri-food supply
	chain.
	\newblock {\em Business Process Management Journal}, 30(1):291--316. https://doi.org/10.1108/BPMJ-01-2023-0035
	
	\bibitem{un2015}
	United Nations (UN) (2015).
	\newblock Transforming our world: The 2030 agenda for sustainable development. https://sdgs.un.org/2030agenda
	
	\bibitem{wan2021novel}
	Wan, S., and Dong, J. (2021).
	\newblock A novel extension of best-worst method with intuitionistic fuzzy
	reference comparisons.
	\newblock {\em IEEE Transactions on Fuzzy Systems}, 30(6):1698--1711. https://doi.org/10.1109/TFUZZ.2021.3064695
	
	\bibitem{wolfert2017big}
	Wolfert, S., Ge, L., Verdouw, C., and Bogaardt, M.~J. (2017).
	\newblock Big data in smart farming--a review.
	\newblock {\em Agricultural systems}, 153:69--80. https://doi.org/10.1016/j.agsy.2017.01.023
	
	\bibitem{wu2023analytical}
	Wu, Q., Liu, X., Zhou, L., Qin, J., and Rezaei, J. (2023).
	\newblock An analytical framework for the best-worst method.
	\newblock {\em Omega}, 123:102974. https://doi.org/10.1016/j.omega.2023.102974
	
	\bibitem{xu2024some}
	Xu, Y., and Wang, D. (2024).
	\newblock Some methods to derive the priority weights from the best--worst
	method matrix and weight efficiency test in view of incomplete pairwise
	comparison matrix.
	\newblock {\em Fuzzy Optimization and Decision Making}, 23(1):31--62. https://doi.org/10.1007/s10700-023-09410-w
	
	\bibitem{youssef2020integrated}
	Youssef, A.~E. (2020).
	\newblock An integrated mcdm approach for cloud service selection based on
	topsis and bwm.
	\newblock {\em IEEE Access}, 8:71851--71865. https://doi.org/10.1109/ACCESS.2020.2987111
	
	\bibitem{zhang2022circular}
	Zhang, Q., Dhir, A., and Kaur, P. (2022).
	\newblock Circular economy and the food sector: A systematic literature review.
	\newblock {\em Sustainable Production and Consumption}, 32:655--668. https://doi.org/10.1016/j.spc.2022.05.010
	
	\bibitem{zhao2018comprehensive}
	Zhao, H., Guo, S., and Zhao, H. (2018).
	\newblock Comprehensive benefit evaluation of eco-industrial parks by employing
	the best-worst method based on circular economy and sustainability.
	\newblock {\em Environment, development and sustainability}, 20(3):1229--1253. https://doi.org/10.1007/s10668-017-9936-6
	
\end{thebibliography}

\end{document}